\documentclass[a4paper, 12pt]{article}

\usepackage[sort&compress]{natbib}
\bibpunct{(}{)}{;}{a}{}{,} 

\usepackage{amsthm, amsmath, amssymb, mathrsfs, multirow, url, subfigure}
\usepackage{graphicx} 
\usepackage{ifthen} 
\usepackage{amsfonts}
\usepackage[usenames]{color}
\usepackage{fullpage}

\theoremstyle{plain} 
\newtheorem{thm}{Theorem}

\theoremstyle{definition}
\newtheorem{definition}{Definition}

\theoremstyle{remark}

\newtheorem{example}{Example}

\newcommand{\prob}{\mathsf{P}}

\newcommand{\pl}{\mathsf{pl}}

\newcommand{\model}{\mathscr{P}}

\newcommand{\pois}{{\sf Pois}}
\newcommand{\bin}{{\sf Bin}}
\newcommand{\unif}{{\sf Unif}}
\newcommand{\nm}{{\sf N}}

\newcommand{\gam}{{\sf Gamma}}

\newcommand{\chisq}{{\sf ChiSq}}

\newcommand{\A}{\mathcal{A}} 
\newcommand{\B}{\mathcal{B}} 
\newcommand{\RR}{\mathbb{R}}

\newcommand{\YY}{\mathbb{Y}}
\newcommand{\UU}{\mathbb{U}}

\newcommand{\M}{\mathscr{M}}

\renewcommand{\S}{\mathcal{S}}

\newcommand{\T}{\mathcal{T}}

\newcommand{\iid}{\overset{\text{\tiny iid}}{\,\sim\,}}

\title{False confidence, non-additive beliefs, and valid statistical inference}
\author{Ryan Martin\footnote{Department of Statistics, North Carolina State University. Email: {\tt rgmarti3@ncsu.edu}}}
\date{\today}

\begin{document}

\maketitle 

\begin{abstract} 
Statistics has made tremendous advances since the times of Fisher, Neyman, Jeffreys, and others, but the fundamental and practically relevant questions about probability and inference that puzzled our founding fathers remain unanswered.  To bridge this gap, I propose to look beyond the two dominating schools of thought and ask the following three questions: what do scientists need out of statistics, do the existing frameworks meet these needs, and, if not, how to fill the void?  To the first question, I contend that scientists seek to convert their data, posited statistical model, etc., into calibrated degrees of belief about quantities of interest.  To the second question, I argue that any framework that returns additive beliefs, i.e., probabilities, necessarily suffers from {\em false confidence}---certain false hypotheses tend to be assigned high probability---and, therefore, risks systematic bias.  This reveals the fundamental importance of {\em non-additive beliefs} in the context of statistical inference.  But non-additivity alone is not enough so, to the third question, I offer a sufficient condition, called {\em validity}, for avoiding false confidence, and present a framework, based on random sets and belief functions, that provably meets this condition.  Finally, I discuss characterizations of p-values and confidence intervals in terms of valid non-additive beliefs, which imply that users of these classical procedures are already following the proposed framework without knowing it.  

\smallskip

\end{abstract}


\section{Introduction}
\label{S:intro}

Statistics is concerned with the collection and analysis of data for the purpose of evaluating existing theories and creating new ones, i.e., learning about the world through observation.  This is central to the scientific method and, therefore, a solid foundation of statistical inference---a set of core principles upon which an effective statistical analysis is built---is essential to the advancement of science.  This foundation has been severely shaken in recent years, however, starting perhaps with Ioannidis's provocative claim that {\em most published research findings are false} \citep{why.most.false}.  Science is a ``culture of doubt'' \citep{feynman.doubt}, so Ioannidis's skepticism is healthy.  Unfortunately, the support for his claim has continued to mount---see \citet{nuzzo2014} and the reports in {\em Nature}\footnote{\url{http://www.nature.com/nature/focus/reproducibility/index.html}}---and the gravity of the problem is now becoming clear: large-scale replication projects have found that just $37/97=39$\% of findings in psychology \citep{open.science.2015}, $11/18 = 61$\% in experimental economics \citep{Camerer2016}, and $13/21=62$\% in social science \citep{camerer2018} could be replicated.  Rightfully so, this plight is now being called a {\em replication crisis}. 

There are many factors that contribute to the replication crisis, some are statistics-related and others are not.\footnote{Sociological factors, such as the ``publish or perish'' culture of academia, are major contributors to the replication crisis; see \citet{crane.martin.trust, crane.martin.mission, crane.martin.quillette}.}  While statisticians can't fully resolve the replication crisis, it's important that we acknowledge our contribution to the problem and do our part to help.  There are many ways to help but, since the crisis reflects poorly on the core of our subject, I believe that foundations require our direct attention.

Sadly, it's necessary to emphasize that discussions of foundations are practically relevant, so let me kick off this discussion in the same way as \citet{barnard.book.1985}:
\begin{quote}
{\em I shall be concerned with the foundations of the subject. But in case it should be thought that this means I am not here strongly concerned with practical applications, let me say that confusion about the foundations of the subject is responsible, in my opinion, for much of the misuse of statistics that one meets in fields of application\ldots}
\end{quote}
Barnard's need to establish the practical relevance of foundational discussions stems from their tendency to turn abstract and philosophical, blurring their connection to everyday statistical practice.  Both abstraction and philosophy are necessary but, to be effective, they must be applied at the right level and in the right direction.  The two dominant schools of thought, namely, {\em frequentist} and {\em Bayesian}, have co-existed for a century despite many philosophical arguments for one approach over the other.  The reason these arguments have failed to take hold is that their practical implementation generally doesn't align with the intuition and experience of data analysts.  An important example is {\em Birnbaum's theorem}\footnote{There are two generally accepted notions in statistics: the {\em sufficiency principle}, which says that inferences should depend only on the value of the minimal sufficient statistic, and the {\em conditionality principle}, which says, roughly, inference should be conditioned on the observed values of ancillary statistics.  A more controversial notion is the {\em likelihood principle}, which says that inferences should only depend on the shape of the likelihood function.  Birnbaum's theorem states that the sufficiency and conditionality principles imply the likelihood principle.  This is problematic because really only a Bayesian approach satisfies the likelihood principle.  Therefore, if I use, say, a $t$-test, which violates the likelihood principle, then Birnbaum's theorem says I'm making an illogical step somewhere because I must have violated the well-accepted sufficiency or conditionality principles.} \citep[e.g.,][]{birnbaum1962, bergerwolpert1984}, which implies that the use of, say, the $z$- or $t$-test about an unknown normal mean is wrong or lacks foundational support despite its non-controversial theoretical properties and many years of successful applications.  Needless to say, practitioners doubt the relevance of such a result, so its impact has been rather limited.\footnote{In addition to the doubt about the relevance of Birnbaum's theorem, there has been doubt about its factuality.  Indeed, \citet{evans2013} and \citet{mayo2014} have argued that the claim is wrong or at least that its implications have been exaggerated.}  

This lack of a clear and agreed-upon foundation or set of core principles is unsettling.  On the one hand, this gives non-statisticians the impression that we don't really know what we're doing; indeed, as \citet{fraser2011.rejoinder} paraphrases from an astute outsider, ``Why don't [statisticians] put their house in order?''  On the other hand, statisticians have mostly given up on trying to sort this out, e.g., 
\begin{quote}
{\em Basic Principle 0: Do not trust any principle}~\citep{lecam.mle.1990},
\end{quote}
and 
\begin{quote}
{\em If a statistical analysis is clearly shown to be effective [...]~it gains nothing from being described as principled}~\citep{speed.principles}.  
\end{quote}
To bridge this gap between statisticians and non-statisticians, I believe that we need to ask different questions.  In particular, statisticians should not be asking {\em frequentist or Bayesian?}~but, rather, {\em how best to meet the needs of science?}  


In this direction, and in the spirit of overcoming the replication crisis described above, I like what \citet{reid.cox.2014} have suggested:
\begin{quote}
{\em \ldots for any analysis that hopes to shed light on the structure of the problem, modeling and calibrated inferences \ldots seem essential.}
\end{quote}
My focus is the inference aspect.  This is not because modeling is easy or of lesser importance, but because the quantity about which inferences are sought is often defined only relative to the posited model.  Since the inference problem needn't exist without a model, I will assume that a model is already in place, and focus on the ``calibrated inferences seem essential'' part.  To make this suggested principle operational, a precise understanding of what each word means is necessary.  
\begin{itemize}
\item {\em Inference}.  There is a strong desire to quantify uncertainties via probabilities or, more generally, via degrees of belief, and to make inferences based on these, so this will be my approach here; see Definition~\ref{def:inference}.  Moreover, I'll argue in Section~\ref{S:converse} that my choice to formulate inferences based on degrees of belief is more inclusive than exclusive, i.e., the classical approaches based on p-values and confidence intervals are actually part of the proposed framework. 
\vspace{-2mm}
\item {\em Calibrated}. In metrology, calibration refers to the comparison of measurements to a specified scale.  In the present context, it's the data analyst's degrees of belief that require a scale for interpretation, e.g., what belief values are sufficiently ``large'' or ``small'' that justified conclusions can be drawn.  Another side of calibration, which appears in a follow-up quote from Reid and Cox (see Section~\ref{S:valid} below), is that the scale on which beliefs are interpreted should prevent ``systematically misleading conclusions,'' i.e., provide control on the frequency of errors.  
\vspace{-2mm}
\item {\em Essential}. This word in the trio I take very seriously, for reasons that I will explain.  Academic statisticians are encouraged/incentivized to focus on ``methodology,'' the development of methods and software to be used (hopefully widely) by practitioners and scientists.  Therefore, the system is set up so that those with the expertise to develop methods are involved in data analyses only at the inference step, and even this involvement is indirect through the practitioner's use of their method/software.  If this is the extent of statisticians' involvement in most scientific applications, should we be satisfied with methods that provide ``approximate calibration'' for ``most hypotheses''?  I will argue in Section~\ref{S:additive} that there are non-trivial dangers lurking, creating risk of systematic bias, so here I'll interpret ``essential'' literally, that calibrated beliefs are absolutely necessary.
\end{itemize}

The main goal of this paper is two-fold.  First, I explore the statistical implications of the suggested ``calibrated inferences seem essential'' principle.  After some setup and notation, in Section~\ref{S:setup}, I introduce the object of study, namely, an {\em inferential model}, which is a map from data, statistical model, etc., to a function that returns the data analyst's degrees of belief in hypotheses of interest.  These inferential models can take many forms, but the most familiar are those that return {\em additive} degrees of belief, i.e., probabilities.  In Section~\ref{S:additive}, I argue that additive inferential models are always at risk of miscalibration and systematic bias and, therefore, apparently don't align with the ``calibrated inferences seems essential'' standard.  This suggests that there is something fundamental about non-additive beliefs in this context.  The presence of non-additivity can be felt in the logic of classical statistical inference---e.g., rejecting the null hypothesis does not mean accepting the alternative---but, to my knowledge, the extent of non-additivity's importance has yet to be fully elucidated.  


This leads to the second main goal of the paper, moving beyond those familiar additive inferential models to something that does meet the specified standard.  Unfortunately, non-additivity alone is not enough.  Additional constraints on the inferential model output are needed, beyond non-additivity, and Section~\ref{S:valid} presents what I call the {\em validity condition}, a mathematical formulation of the requirements contained in the phrase ``calibrated inferences seem essential.''  After some important remarks in Section~\ref{S:remarks1}, and a review of non-additive beliefs in Section~\ref{S:nonadditive}, I describe the framework from \citet{imbook} that is able to meet the validity condition.  What distinguishes this approach from others is the incorporation of a suitable, user-specified random set, which results in a provably valid inferential model that returns a genuinely non-additive belief function.  General details about the construction, the validity theorem, and a number of examples are presented in Section~\ref{S:valid.im}.  There I also describe some characterization results revealing that users of those classical statistical methods, such as p-values and confidence intervals, have actually unknowingly been working in my non-additive belief framework all along.  After another collection of remarks in Section~\ref{S:remarks2}, I conclude in Section~\ref{S:discuss} with some perspectives, open problems, etc.

\section{Inferential models}
\label{S:setup}

\subsection{Setup, notation, and objectives}

To start, I assume that there is a statistical model $\prob_{Y|\theta}$, indexed by a parameter $\theta \in \Theta$, for the observable data $Y \in \YY$.  
Here, as is common in discussions of inference, I will assume that both the data and model are given, but they both deserve a few comments; see also \citet{crane.martin.needs}.  

\begin{itemize}

\item Data is what the statistician gets to see; in fact, it's the only thing in a data analysis that's really known.  Data must certainly be relevant to the scientific problem, and the more informative the better.  Technical details about the collection of (relevant and informative) data are presented in \citet{hedayat.sinha.1991} and \citet{hinkelmann.kempthorne.2008}.  Here I use the adjective ``observable'' to distinguish $Y$ from other random variables that are ``unobservable;'' see Section~\ref{S:valid.im}.  

\vspace{-2mm}

\item \citet[][Ch.~1]{lehmann.casella.1998} define a model as a family of probability measures on $\YY$, indexed by a parameter $\theta \in \Theta$:
\begin{equation}
\label{eq:model}
\model = \{\prob_{Y|\theta}: \theta \in \Theta\}.
\end{equation}
(Here $\prob_{Y|\theta}$ is a probability distribution for $Y \in \YY$, a random variable, random vector, or whatever, depending on $\theta$; it need not be interpreted as a conditional distribution of $Y$, given $\theta$, derived from a joint distribution for the pair.) Guidelines on the development of sound models are given, e.g., in \citet{cox.hinkley.book}, \citet{box1980}, and \citet{mccullagh.nelder.1983}; \citet{mccullagh2002} gives a more technical discussion.  Here I will assume that the statistical model \eqref{eq:model} is derived from some knowledge about the data-generating process and/or from an exploratory data analysis.  This model serves two important purposes: first, it establishes a formal connection between the data and the real-world problem; and, second, the relevant scientific questions can be formulated precisely in terms of the model parameter, $\theta$.  While it's not a focus of the present paper, I should say that the development of sound models is far from routine, especially in modern, complex problems such as network analysis \citep[e.g.,][]{crane.book}.  
\end{itemize}

Some might say that the model is incomplete in the sense that it's missing an uncertainty assessment about the parameter $\theta$ in the form of a prior; see, for example, \citet{definetti.book.1972}, \citet{savage1972}, \citet{gelman.bda.2004}, and \citet{kadane.book}.  If a real subjective prior is available, then the statistical inference problem is ``simply an application of probability'' \citep[][p.~xxv]{kadane.book}, in particular, Bayes's theorem \citep[e.g.,][]{efron2013.science} yields the conditional distribution of $\theta$, given $Y=y$, which can be used for inference.  If genuine prior information is available, then surely one should use it.  However, a substantial part of the literature adopts a view that no real prior information is available; even in the Bayesian literature, the trend is to work with ``default'' or ``non-informative'' priors that let the data speak for itself.  One reason for this interest in prior-free scenarios is, as Brad Efron said during a presentation at the 2016 Joint Statistical Meetings in Chicago, that ``scientists like to work on new problems.''  In any case, my choice to separate the prior and model is to explore what can be done in cases where a fully satisfactory prior distribution is not available.  Incorporating prior information into the framework I'm describing here is possible and I'll have more to say on this in Section~\ref{S:discuss}.   

Given a sound statistical model and relevant data, the goal is to make inferences about the unknown $\theta$.  Classically, these inferences come in the form of point or interval estimates, hypothesis tests, etc.; prediction could also be included in this list of tasks \citep[e.g.,][]{impred}, but my focus here is on inference about fixed-but-unknown quantities.  A selling point of Bayesian and other distribution-based approaches to statistical inference is that it's relatively straightforward to read off, for example, point and interval estimates, from the posterior distribution.  So, in the same spirit, I will take this ``posterior distribution'' or, more generally, the data analyst's degrees of belief (based on data and posited statistical model) as the primitive, with estimates and tests being derivatives.  The next subsection makes this more precise.

\subsection{Definition}

In the Bayesian literature \citep[e.g.,][]{kadane.book, lindley.book}, those elements discussed in the previous section are converted to a probability distribution on $\Theta$ that encodes the data analyst's degrees of belief, from which inferences will be drawn.  In other words, the Bayesian procedure, or {\em model}, for inference is one in which the various inputs---relevant data, sound statistical model, etc.---are mapped to probabilities concerning the unknown quantities.  This idea can be generalized, as in the following definition, allowing for degrees of belief output that need not be a probability distribution.   

\begin{definition}
\label{def:inference}
Fix the sample space $\YY$, parameter space $\Theta$, and let $\model$ be a statistical model as in \eqref{eq:model} that connects the two.  Then an {\em inferential model} is a map from $(\model, y, \ldots)$ to a function $b_y: 2^\Theta \to [0,1]$ where, for each hypothesis $A \subseteq \Theta$, $b_y(A)$ represents the data analyst's degree of belief in the truthfulness of the assertion ``$\theta \in A$,''  and ``$\ldots$'' denotes extra inputs that may vary from one situation to another.  
\end{definition}

This definition formalizes the idea that inference is based on a process of converting observations into degrees of belief about $\theta$.  That is, an inferential model is simply a rule or procedure that describes how the observed data, posited statistical model, and perhaps other ingredients (e.g., a prior distribution) are processed to quantify uncertainty for the purpose of making inference.  In particular, if the relevant scientific question is encoded in terms of some feature $\phi(\theta)$ of the model parameter $\theta$, then the goal is to evaluate degrees of belief $b_y(A)$ for hypotheses of the form $A = \{\theta: \phi(\theta) \in B\}$ for suitable $B \subseteq \phi(\Theta)$.   

Note that the inferential model and the corresponding degrees of belief output, $b_y$, depend directly on $\model$.  A reviewer interpreted the form of this dependence as suggesting that the inferential model is just a wrapper around the inputs $(\model, y)$ as opposed to all the pieces contributing equally.  I agree with this interpretation, but I don't see this as unexpected or problematic.  In applications, after collecting data and formulating a model, the investigator still has options for how to analyze the data, e.g., construct a Bayesian posterior distribution, carry out a likelihood ratio test, etc.  And, again, that there is an order---first model, then inference---is necessary because the precise formulation of the relevant scientific questions is in terms of the model parameters, so the inference problem generally isn't even well-defined before the model is specified. 

This inferential model notion captures the familiar Bayesian approach, where the output $b_y$ is simply the posterior distribution for $\theta$, given data $y$, under the assumed model, with a user-specified prior distribution as part of ``...''  It also captures fiducial inference \citep{fisher1973, zabell1992, barnard1995}, {structural inference} \citep{fraser1968}, {generalized inference} \citep{chiang2001, weerahandi1993}, {generalized fiducial inference} \citep{hannig2009, hannig2012, hannig.lee.2009, hannig.iyer.patterson.2006, lai.hannig.lee.2013, hannig.review}, {confidence distributions} \citep{schweder.hjort.2002, schweder.hjort.book, xie.singh.2012}, and Bayesian inference with data-dependent priors \citep{fraser.reid.marras.yi.2010, martin.walker.eb, martin.walker.deb}.  The key feature these all share is that their output, $b_y$, is a probability measure.  But Definition~\ref{def:inference} does not require $b_y$ to be a probability, so it also covers situations where $b_y$ is something more general, such as a capacity, belief function, etc., as in Section~\ref{S:nonadditive}.  In these latter approaches, $b_y$ is not additive, so it is possible that {\em both} $b_y(A)$ and $b_y(A^c)$ are small if data $y$ is not especially informative about $A$.  It turns out that additive versus non-additive has important consequences; see Section~\ref{S:additive}.  

Whatever the specific mathematical form of $b_y$, it can be used in a natural and straightforward way.  That is, large values of $b_y(A)$ and $b_y(A^c)$ suggest that data strongly supports the truthfulness and falsity, respectively, of the claim ``$\theta \in A$.''\footnote{Intermediate values of $b_y(A)$, such as 0.4, are more difficult to interpret, but the same is true for probabilities; see Section~\ref{SS:beyond}.}  See \citet{dempster.copss} for a discussion of this general approach to inference through degrees of belief.  Moreover, the data analyst can, if desired, produce decision procedures based on the inferential output, e.g., reject a null hypothesis ``$\theta \in A$'' if and only if $b_y(A^c)$ is sufficiently large.  

Lastly, it should be clear that inferential models are not unique---my degrees of belief might be different from yours.  That inference requires individual judgment is not a shortcoming, it is what makes statistics an interesting subject and what makes statisticians' expertise valuable.  Depending on the application, however, I may want my degrees of belief to be meaningful to others, not just to myself, which puts certain constraints on my inferential model.  As \citet[][p.~xvi]{evans2015book} writes, ``subjectivity can never be avoided but its effects can be\ldots controlled.''  This control comes from external or, in some sense, ``non-subjective'' considerations, and the remainder of the paper focuses on these.

\section{Additive inferential models}
\label{S:additive}

\subsection{Satellite conjunction analysis}
\label{SS:conjunction}

Satellites orbiting Earth impact our everyday lives.  For example, those like me with a poor sense of direction rely on GPS navigation applications that use positioning information from satellites to successfully reach our destinations.  But it is not enough just to get these satellites into orbit, they must be constantly monitored to prevent collision with one another or with other kinds of debris.  The challenge is that collision events cannot be determined with certainty, so conjunction analysts must make important decisions, e.g., whether evasive maneuvers are necessary, in the face of uncertainty.  This practically-relevant application will serve as motivation for the forthcoming general discussion of additive versus non-additive beliefs.  

Conjunction analysts to use probabilities for uncertainty quantification, and their approach is as follows.  Consider a pair of orbiting satellites whose positions and velocities are points in three-dimensional space.  The standard conjunction analysis starts with an approximation to reduce it to a two-dimensional problem.  That is, a plane of closest approach is determined and then it is assumed that the satellites are each approaching this plane from an orthogonal direction, so only the projections onto this plane are relevant.  Let $\theta=(\theta_1,\theta_2)^\top$ represent the difference between the true but unknown satellite positions on this plane, so that small $\|\theta\|$ indicates the satellites are close, hence a collision.  Let $Y=(Y_1,Y_2)^\top$ denote the estimate of this difference in satellite positions based on noisy data.  A simple statistical model to describe this scenario is $Y \sim \nm_2(\theta, \Sigma)$. The normality assumption is often justified by invoking the central limit theorem, since the estimate $Y$ is the result of aggregating many position and velocity measurements.  Here, as is common in the conjunction analysis literature, I will assume that the covariance matrix $\Sigma$ can be accurately estimated so it will be taken as fixed.   

The next step in the analysis, as reviewed in \citet{balch.martin.ferson.2017}, is to first construct a corresponding posterior distribution for $\theta$.  Since $\theta$ is a location parameter, one might consider an inferential model that describes degrees of belief based on 
\[ (\theta \mid y) \sim \Pi_y := \nm_2(y, \Sigma). \]
This is a flat-prior Bayes posterior and a fiducial/confidence distribution.  Using information on the satellites' sizes, a collision radius, $t$, can be derived, so that a collision event is defined as ``$\|\theta\| \leq t$.''  The analysis proceeds by calculating the probability of collision (or non-collision) based on the posterior distribution above.  While experienced statisticians may have reasons for concern (see below), this approach is not obviously unreasonable.  

\begin{figure}[t]
\begin{center}
\subfigure[$\|y\|^2 = 0.5$ vs $t=1$]{\scalebox{0.6}{\includegraphics{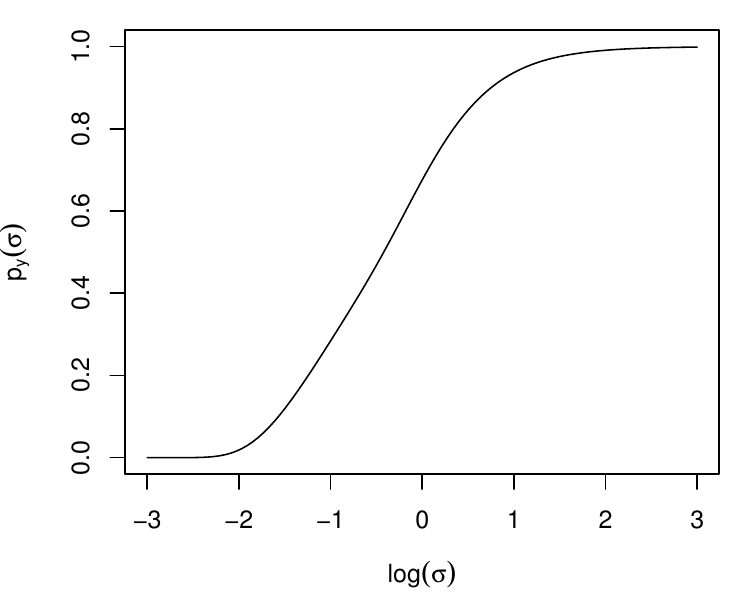}}}
\subfigure[$\|y\|^2 = 0.9$ vs $t=1$]{\scalebox{0.6}{\includegraphics{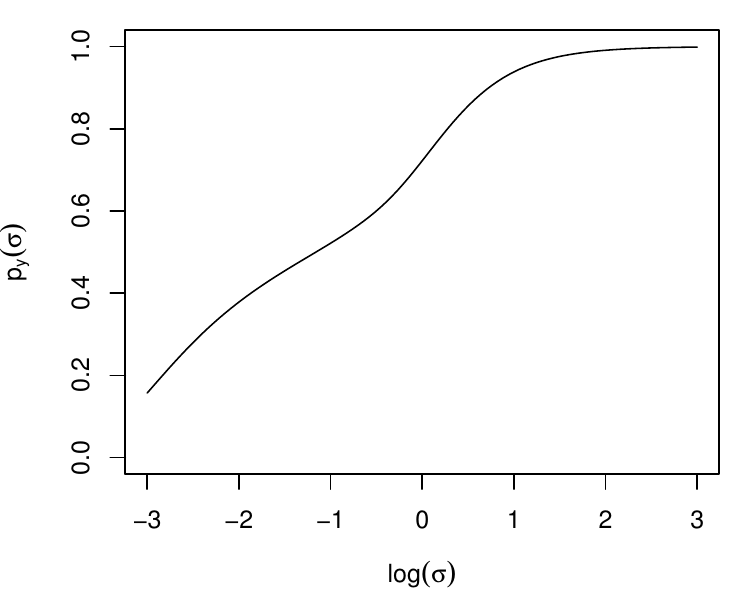}}}
\subfigure[$\|y\|^2 = 1.1$ vs $t=1$]{\scalebox{0.6}{\includegraphics{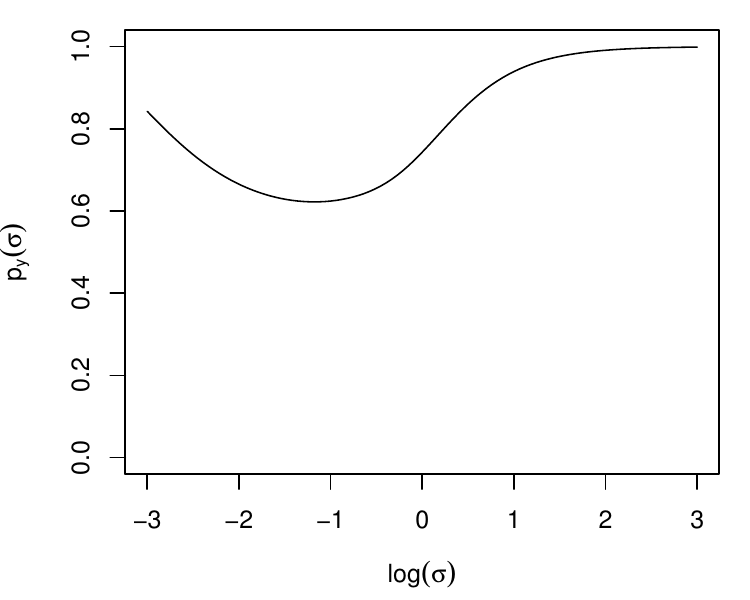}}}
\subfigure[$\|y\|^2 = 1.5$ vs $t=1$]{\scalebox{0.6}{\includegraphics{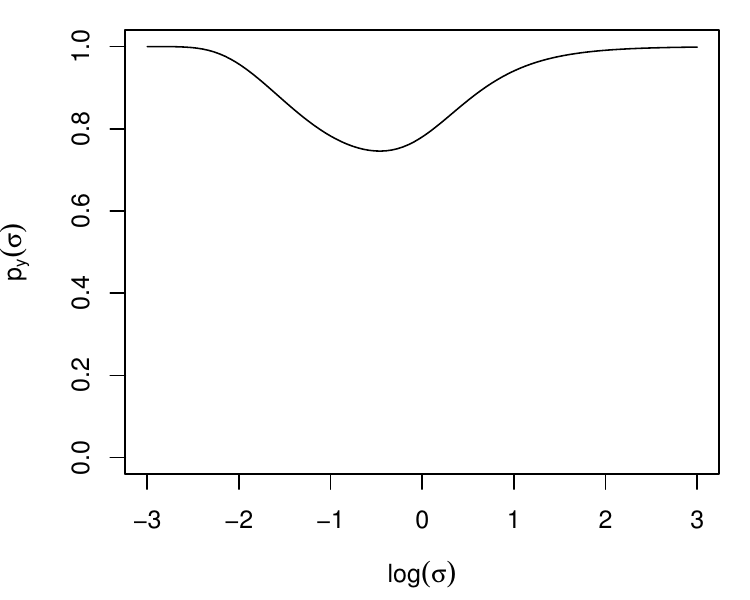}}}
\end{center}
\caption{Plots of the non-collision probability in \eqref{eq:ncp} as a function of $\sigma$ for several configurations of $y$.  Here the collision threshold is $t=1$.}
\label{fig:collision}
\end{figure}

However, an interesting and paradoxical phenomenon has been observed in the conjunction analysis literature, called {\em probability dilution}.  The basic idea is that lower data quality (i.e., more noise in the position and velocity measurements) implies lower estimated collision probability.  This is paradoxical because lower quality data clearly means more uncertainty but the probability calculation gives a different conclusion, namely, that there is more certainty in the satellite's safety.  Of course, if conjunction analysts take a low collision probability at face value, then they are likely to conclude that the satellite is safe when, in fact, it may not be.  To see this clearly, consider the simple scenario where $\Sigma = \sigma^2 I_2$.  Then the non-collision probability, according to the aforementioned posterior distribution, is a simple non-central chi-square calculation:
\begin{equation}
\label{eq:ncp}
\Pi_y(\|\theta\| > t) = 1 - {\tt pchisq}(t^2/\sigma^2, {\tt df}=2, {\tt ncp}=\|y\|^2/\sigma^2). 
\end{equation}
Figure~\ref{fig:collision} plots this non-collision probability as a function of $\sigma$.  When $\sigma$ is small, there is no problem; but when $\sigma$ is large, the non-collision probabilities are large, {\em no matter what the data $y$ says}.  That is, even if the data suggests a collision, like in the $\|y\|^2=0.5$ case, with enough noise, the posterior probabilities indicate that the satellite is safe.  This turns out to be a special case of the {\em false confidence} phenomenon discussed in Sections~\ref{SS:cv}--\ref{SS:fct}.  

Statisticians will recognize this problem---inference on the length of a multivariate normal mean vector---from \citet{stein1959}, so the paradoxical conclusions might be expected.  But this normal means problem is not the only one where issues like this arise, so a better understanding of the causes and effects is valuable.  Besides, the insights gleaned from Stein's investigation of this deceptively simple problem, in particular, the notions of {\em shrinkage} and {\em borrowing of strength}, are now central to modern high-dimensional data analysis, so it is worth digging into these peculiar examples further.

\subsection{Normal coefficient of variation}
\label{SS:cv}

For data $Y=(Y_1,\ldots,Y_n)$, consider a statistical model that posits the components of $Y$ to be independent and identically distributed (iid), with $Y_i \sim \nm(\mu, \sigma^2)$, $i=1,\ldots,n$.  Here $\theta=(\mu, \sigma^2)$ is the unknown parameter, but primary interest is in $\phi = \sigma / \mu$, a quantity called the {\em coefficient of variation}.  One possible approach is to construct a Bayesian inferential model where, given a prior distribution for $\theta$, the corresponding posterior for $\theta$ is obtained via Bayes's formula, and then a marginal posterior for $\phi$ is derived from that of $\theta$ using the ordinary probability calculus.  

Since interest is in $\phi$, I will adopt the reference prior for $\theta$, now expressed as $\theta=(\phi, \sigma)$, presented in \citet[][Example~7]{berger.liseo.wolpert.1999}, designed specifically for cases where inference on $\phi$ is the goal.  Using the standard default conditional prior $\pi(\sigma \mid \phi) \propto \sigma^{-1}$, for $\sigma$, given $\phi$, and the reference prior for $\phi$, 
\[ \pi(\phi) \propto \{|\phi| (\phi^2 + \tfrac12)^{1/2} \}^{-1}, \]
given in their Equation~(38), one obtains the proper marginal posterior for $\phi$:
\begin{equation}
\label{eq:cv.post}
\pi_y(\phi) \propto \pi(\phi) e^{-\frac{n}{2\phi^2}(1 - \frac{\bar y^2}{D^2})} \int_0^\infty z^{n-1} e^{-\frac{nD^2}{2}(z - \frac{\bar y}{D^2\phi})^2} \, dz, 
\end{equation}
where $\bar y = n^{-1} \sum_{i=1}^n y_i$ and $D^2 = n^{-1} \sum_{i=1}^n y_i^2$.  This agrees with their Equation~(39).  The integral can be evaluated numerically, so posterior computations via Markov chain Monte Carlo, e.g., the Metropolis--Hastings algorithm, are relatively straightforward.  

\begin{figure}
\begin{center}
\scalebox{0.7}{\includegraphics{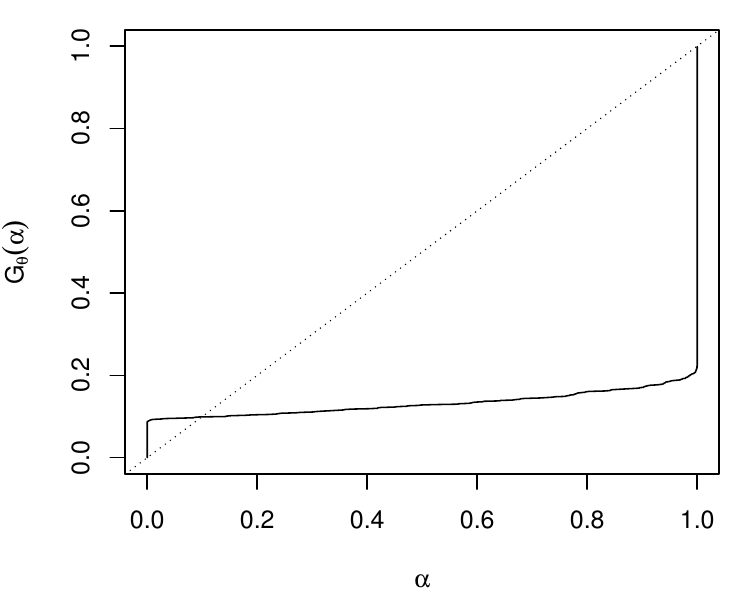}}
\end{center}
\caption{Plot of the estimated distribution function $G_\theta(\alpha) = \prob_{Y|\theta}\{\Pi_Y(A) \leq \alpha\}$ of the Bayesian posterior probability, $\Pi_Y(A)$, where $A=(-\infty,9]$ is a false hypothesis about the coefficient $\phi=\sigma/\mu$ when $\mu=0.1$ and $\sigma=1$.}
\label{fig:cv1}
\end{figure}

However, like in the satellite conjunction analysis example from the previous section, there is no guarantee that this posterior probability distribution for the interest parameter provides reliable inference.  Here I will explore this question using a simulation study.  Let $Y$ consist of $n=10$ iid samples from $\nm(\mu, \sigma^2)$, with $\mu=0.1$ and $\sigma=1$; then the true value of $\phi$ is $10$.  Next, consider a hypothesis about $\phi$, namely, $A=(-\infty, 9]$, which happens to be false.  Let $\Pi_y(A)$ denote the  probability of this hypothesis with respect to the posterior distribution in \eqref{eq:cv.post}.  Then $\Pi_Y(A)$ is a random variable, as a function of $Y$, and I am interested in its distribution function, $G_\theta(\alpha) = \prob_{Y|\theta}\{\Pi_Y(A) \leq \alpha\}$.  Figure~\ref{fig:cv1} shows the estimated distribution function, based on 1000 data sets.  The most striking feature of this plot is the high concentration of the distribution's mass around value 1.  This is problematic because, remember, this is the posterior probability of a {\em false} hypothesis about $\phi$.  In those instances when $\Pi_y(A) \approx 1$ and, therefore, $\Pi_y(A^c) \approx 0$, the only reasonable conclusion would be that the data support hypothesis $A$.  Since these problematic instances are apparently frequent, inferences drawn based on the posterior distribution described above would be systematically misleading.  

It is worth mentioning that there is nothing special about the particular settings in this illustration.  For any values of $(n, \mu, \sigma)$, there exists a hypothesis $A$, perhaps of similar form to the one above, such that, even if $A$ is false, the posterior probability assigned to it will tend to be large; see Section~\ref{SS:fct}.  

Some statisticians might also recognize this normal coefficient of variation problem as one in the collection of challenging inference problems identified by \citet{gleser.hwang.1987} and \citet{dufour1997} and discussed in \citet{bertanha.moreira.impossible}.  One could say that this and others in these classes are not a practical problems, and that might be true.  But these ``impossible'' inference problems are precisely the ones we should use to test our understanding, in the same way that we test our understanding of measure theory by constructing and examining non-measurable sets, etc.

\subsection{False confidence theorem}
\label{SS:fct}

The previous two subsections showed that there are instances where inferences based on a Bayesian inferential model would be systematically misleading.  It turns out that these issues have nothing to do with the Bayesian features of those solutions.  Instead, false confidence is a consequence of using {\em additive} degrees of beliefs, i.e., probabilities, to describe uncertainty.  The following theorem makes this precise.  

\begin{thm}[\citealt{balch.martin.ferson.2017}]
\label{thm:fct}
Consider probability measure $\Pi_y$ with a density function that is bounded and continuous for each $y$.  Then for any $\theta \in \Theta$, any $\alpha \in (0,1)$, and any $p \in (0,1)$, there exists a set $A \subset \Theta$ such that 
\[ 
A \not\ni \theta \quad \text{and} \quad \prob_{Y|\theta}\{\Pi_Y(A) \geq 1 - \alpha \}  \geq p. \]
\end{thm}

In words, any inferential model with additive degrees of belief---including Bayes, fiducial, etc.---will be at risk of false confidence.  That is, there are false hypotheses to which $\Pi_y$ will tend to assign high probability, hence a risk of systematic biase.  It is important to emphasize that the theorem only says {\em there exists} hypotheses afflicted by false confidence.  Certainly, some of these hypotheses are trivial, e.g., complements of very small subsets of $\Theta$, but the examples in Sections~\ref{SS:conjunction}--\ref{SS:cv} reveal that non-trivial hypotheses can be affected too; see, also, \citet{williams.fct.2018}.  In my experience, the non-trivial yet problematic hypotheses have unusual shapes: the complement of a disc in the plane and the complement of a cone in the upper half-plane in Section~\ref{SS:conjunction} and \ref{SS:cv}, respectively.  These unusual shapes most commonly arise when hypotheses concern a non-linear function of the full parameter.  Similar concerns about the performance of Bayesian posterior inference for non-linear functions of $\theta$ were raised in \citet{fraser2011} and \citet{fraser.etal.2016}, but Theorem~\ref{thm:fct} expands on their conclusions in several directions, in particular, it applies to any additive inferential model, not just Bayesian, and to more general hypotheses.  Moreover, it helps to explain why, as advised by \citet[][p.~66]{cox2006} and others, it is dangerous to manipulate, e.g., marginalize, fiducial and confidence distributions according to the usual probability calculus.  

When genuine prior information is available in the form of a proper prior distribution for $\theta$, so that the behavior of $\Pi_Y(A)$ would be assessed with respect to the marginal distribution of $Y$ derived from the joint distribution of $(Y,\theta)$, then false confidence is also of no concern.  For one thing, there are no {\em false hypotheses} in such a context since $\theta$ is not an unknown parameter in the marginal distribution of $Y$.  Also, there are certain automatic controls on the distribution of $\Pi_Y(A)$ in such settings, but these bounds are generally only as meaningful as the postulated prior; see Section~\ref{SS:bayes.valid}.  

There is also a connection between the false confidence theorem and the non-existence of a truly ``non-informative'' prior distribution established in Proposition~2.1 of \citet{imbook}.  That is, if every prior probability distribution is ``informative'' simply by virtue of its additivity, not because it's based on real prior information, then there will be false hypotheses towards which the posterior distribution is pushed.  

In certain communities, additive probabilities are called {\em precise} in the sense that the postulated statistical model fully captures all the existing uncertainty.\footnote{The {\em 2019 International Symposium on Imprecise Probabilities}, \url{http://isipta2019.ugent.be}, has a theme, related to this point, that I like: ``there's more to uncertainty than probabilities.''}  Such an assumption is always difficult to justify but, for sure, said justification would require that the prior distribution be informative in some sense.  With the commonly used default prior distributions, any claims of ``precision'' would be questionable:
\begin{quote}
\label{quote:fraser.real}
{\em [Bayes's formula] does not create real probabilities from hypothetical probabilities...} \citep[][p.~249]{fraser.copss}
\end{quote}
There are so many hypotheses in the $\sigma$-algebra upon which the probability is defined and, given the necessarily limited informativeness of the data, it is unrealistic to expect that a posterior distribution can assign meaningful or otherwise reliable probabilities to every one.  Therefore, there must be hypotheses about which the posterior probabilities should not be trusted, and this is what the false confidence theorem says.   

On its own, that there exists hypotheses afflicted by false confidence may not be a serious concern.  The problem is that one generally does not know which hypotheses are afflicted and which are not.  Arguably the strongest selling point of Bayesian and other approaches that produce probabilities as output is that they provide a means to answer {\em virtually any question} about the unknowns, specifically, that the full posterior can be marginalized to any quantity of interest to make justifiable inferences.  The false confidence theorem says that this understanding is wrong, that is, there are features of $\theta$ for which certain questions cannot be answered reliably or justifiably; see, also, \citet{grunwald.safe}.  Whether those afflicted hypotheses are of any practical interest depends on the context, but not knowing which ones they are puts the data analyst at risk.  Short of identifying those afflicted hypotheses in each individual situation to mitigate risk, the only option is to avoid false confidence altogether, which, according to Theorem~\ref{thm:fct}, requires breaking additivity.  Not every non-additive belief will be free of false confidence, however, so some additional constraints are needed, as discussed in Section~\ref{S:valid} below.

\section{Validity condition}
\label{S:valid}

If false confidence is undesirable, as I say it is, then an important question is how to avoid it.  Theorem~\ref{thm:fct} says additivity must be abandoned, and I will discuss specific types of non-additive degrees of belief in Section~\ref{S:nonadditive}.  But non-additivity alone is not enough to avoid false confidence so, before getting into those details, here I introduce a condition that will ensure false confidence can be avoided.  The construction of an inferential model that achieves this condition is the topic of Section~\ref{S:valid.im} below.  

The inferential model formulation is one where inferences are drawn according to a data analyst's degree of belief based on data, a statistical model, etc.  In a vacuum, I am free to construct and interpret my degrees of belief in any way that I please.  But if the goal is to communicate my degrees of belief to others, to justify my conclusions, like in a scientific setting, then an agreed-upon scale is needed in order to interpret degrees of belief.  For example, if I conclude that hypothesis $A$ is false because $b_y(A^c) > c$, then why might you find this convincing?  It is not necessary that your degrees of belief agree 100\% with mine, it is enough if I can demonstrate that the event $\{y: b_y(A^c) > c\}$ is rare if $A$ is actually true.  Of course, a more precise connection between the threshold ``$c$'' and the notion of ``rare'' are needed (see Definition~\ref{def:valid}) but the basic idea is familiar and uncontroversial.  Indeed, \citet[][p.~45--46]{cox.hinkley.book} write
\begin{quote}
{\em We should not follow procedures which for some possible parameter values would give, in hypothetical repetitions, misleading conclusions most of the time.}
\end{quote}
Recently, \citet{reid.cox.2014} echoed this point
\begin{quote}
{\em Even if an empirical frequency-based view of probability is not used directly as a basis for inference, it is unacceptable if a procedure\ldots of representing uncertain knowledge would, if used repeatedly, give systematically misleading conclusions}.  
\end{quote}
Both of these remarks mention ``procedures'' which seems to suggest a purely frequentist or behaviorist mode of reasoning based on rules and suggested actions, e.g., reject a hypothesis if and only if such and such happens.  But \citet{fisher1973} and others \citep[e.g.,][]{mayo.cox.2006, mayo.book.1996, mayo.book.2018} argue that there's more going on here, that the frequentist's error rate control leads to an inductive logic for justifying conclusions in individual instances.  Indeed, \citet[][p.~46]{fisher1973} writes
\begin{quote}
{\em ...the feeling induced by a test of significance has an objective basis in that the probability statement on which it is based is a fact communicable to, and verifiable by, other rational minds.  The level of significance in such cases fulfils the conditions of a measure of the rational grounds for the disbelief it engenders.}
\end{quote}
Along similar lines, \citet[][p.~16]{mayo.book.2018} formalizes the notion of an argument from observation, through the absence of error, to conclusion as follows:
\begin{quote}
{\em There is evidence an error is absent to the extent that a procedure with a very high capability of signaling the error, if and only if it is present, nevertheless detects no error.}
\end{quote}
To avoid what Reid and Cox call ``unacceptable,'' and to facilitate the inductive logic of Fisher and Mayo, some kind of calibration is needed, and I suggest the following constraint on ones inferential model.  

\begin{definition}
\label{def:valid}
An inferential model $(\model, y, \ldots) \mapsto b_y$ is {\em valid} if 
\begin{equation}
\label{eq:valid}
\sup_{\theta \not\in A} \prob_{Y|\theta}\{b_Y(A) > 1-\alpha\} \leq \alpha, \quad \forall \; \alpha \in [0,1], \quad \forall \; A \subset \Theta. 
\end{equation}
That is, if $\theta \not\in A$, so that the hypothesis $A$ is false, then $b_Y(A)$, as a function of $Y \sim \prob_{Y|\theta}$, is stochastically no larger than $\unif(0,1)$. 
\end{definition}

The validity property \eqref{eq:valid} implies that, if $A$ is false, then there is a small probability---with respect to the posited statistical model---that $b_Y(A)$ takes large values.  In other words, data are unlikely to provide strong support for false assertions, therefore avoiding Reid and Cox's ``systematically misleading conclusions.'' Since validity prevents the degrees of belief in a false hypothesis from being too large, in a distributional sense, it follows that a valid inferential model is not susceptible to false confidence.  Similar calibration properties are discussed in \citet{balch2012} and \citet{denoeux.li.2018}.  

Of course, since probabilities are susceptible to false confidence, additive inferential models must automatically fail to be valid in the sense of Definition~\ref{def:valid}.  However, this can be seen directly from the fact that \eqref{eq:valid} covers {\em all} assertions, including singletons and their complements.  Indeed, if $A$ is the complement of a singleton and $b_y$ is a probability measure, then usually $b_y(A) = 1$ for all $y$, so \eqref{eq:valid} fails.  But the examples in Section~\ref{S:additive} show that it's not just trivial hypotheses that can be afflicted by false confidence.   

The validity property depends on the assumed statistical model because it imposes a constraint \eqref{eq:valid} on probabilities calculated with respect to $\prob_{Y|\theta}$.  This is not unexpected, since the model defines the parameter about which inferences are desired, so it generally only makes sense to evaluate the quality or reliability of such inferences with respect to the assumed model.  This is not unique to the perspective I'm suggesting here; indeed, the classical size and power of tests and coverage of confidence sets are in the same spirit as \eqref{eq:valid}, i.e., properties of procedures with respect to the assumed model.  

There are two specific questions that I would like to address.  First, {\em why be concerned about all hypotheses?}  Instead, one could ask that the inferential model satisfy a condition like that in \eqref{eq:valid} just for the one or few hypotheses of interest.  There is no problem weakening the condition \eqref{eq:valid} to hold just for a specific collection of hypotheses, see \citet[][Def.~4.2]{imbook}, and indeed this would be satisfactory to the data analyst in a specific application.  But at the meta level where data scientists are developing {\em methods} and corresponding software to used by others across applications, this strong control on performance is necessary.  If I develop a method that automates the inferential model construction, with corresponding software that returns degrees of beliefs for any hypothesis, then unless I identify and warn the user about those problematic hypotheses, I should assume that users will use it for any hypothesis and, therefore, I must protect myself from the associated risk.  For example, if I recommend an inferential model that returns degrees of belief about $\theta=(\mu,\sigma)$ based on samples from $\nm(\mu, \sigma^2)$, and I fail to warn the user about considering hypotheses of the form $A=\{(\mu, \sigma): \sigma \leq c \mu\}$, for given $c \in \RR$, then I share responsibility for any erroneous conclusions based on this.  Consequently, the data scientist developing methods for mass consumption actually has {\em skin in the game} \citep{taleb.skin}, potentially a lot, so to protect him/herself from the associated risk, a strong, uniform control on performance like in \eqref{eq:valid}, for any hypotheses the user might be interested in, would be desired.  

Second, {\em why the specific {\em ``$1-\alpha$''} in \eqref{eq:valid}?}  In other words, \eqref{eq:valid} could be written as 
\[ \sup_{\theta \not\in A} \prob_{Y|\theta}\{b_Y(A) > g_A(\alpha)\} \leq \alpha, \]
for some function $g_A(\alpha)$, so why not allow this more general notion of validity?  Recall that a primary motivation for the validity condition was to have a scale on which the degrees of belief could be interpreted.  To me at least, interpreting probabilities is a challenge because they are {\em relative}.  For example, Nate Silver's {\em FiveThirtyEight} website estimated the (posterior) probability of a Donald Trump victory in the 2016 U.S.~presidential election at 0.28, which was weak enough to predict that Hilary Clinton would win \citep[e.g.,][]{crane.martin.election}, whereas, in, say, the context of variable selection in regression, a (posterior) probability of 0.28 assigned to one of the $2^p$ possible subsets of the $p$ predictor variables would be very strong support for that configuration.  In the same way, the more general definition with $g_A(\alpha)$ in the above display would allow the scale of interpretation to depend on the hypothesis $A$ as well as on other specific features of the problem.  To avoid the confusion of a moving scale of interpretation, an absolute scale is needed, but the specific choice of $\unif(0,1)$ as an absolute scale is not essential---it's just familiar and natural.\footnote{To see what's special about the $\unif(0,1)$ scale, consider the following example.  Suppose you're visiting a place you've never been before and the weatherman on the TV says there's a 10\% chance of rain.  Likely you wouldn't be too concerned about not having an umbrella, but why?  Maybe this weatherman reports {\em only} 0--10\% chance of rain, no matter the circumstances, so that a 10\% chance is actually very high.  That this possibility seems totally absurd suggests that it's ``natural'' to think of probabilities as being calibrated to a $\unif(0,1)$ scale, i.e., that it rains roughly 10\% of the days the weatherman predicts a 10\% chance of rain, and that's why I've imposed this in \eqref{eq:valid}.}  See, also, Section~\ref{SS:bayes.valid}. 

An important consequence of the validity condition is that procedures designed based on the inferential model output will have frequentist guarantees.  To see this clearly, define the dual to the inferential model's degrees of belief as 
\[ p_y(A) = 1 - b_y(A^c), \quad A \subseteq \Theta. \]
Of course, since the focus is on non-additive degrees of belief, $p_y$ is different from $b_y$.  The interpretation of $p_y$ is as a measure of {\em plausibility}, which I discuss further in Section~\ref{S:nonadditive}.  Since \eqref{eq:valid} holds for all $A \subseteq \Theta$, there is an equivalent formulation in terms of $p_y$, i.e., 
\begin{equation}
\label{eq:valid.pl}
\sup_{\theta \in A} \prob_{Y|\theta}\{p_Y(A) \leq \alpha\} \leq \alpha, \quad \forall \; \alpha \in [0,1], \quad \forall \; A \subseteq \Theta. 
\end{equation}
Besides complementing the degrees of belief, $p_y$ provides a simple and intuitive recipe to construct procedures with guaranteed frequentist error rate control.  

\begin{thm}[\citealt{imbasics}]
\label{thm:freq}
Let $(\model, y, \ldots) \mapsto b_y$ be an inferential model that satisfies the validity condition in Definition~\ref{def:valid}. Fix $\alpha \in (0,1)$.  
\begin{enumerate}
\item[{\em (a)}] Consider a testing problem with $H_0: \theta \in A$ versus $H_1: \theta \not\in A$, where $A$ is any subset of $\Theta$.  Then the test, $T_y$, that rejects $H_0$ if and only if $p_y(A) \leq \alpha$ has frequentist Type~I error probability upper-bounded by $\alpha$.  That is, 
\[ \sup_{\theta \in A} \prob_{Y|\theta}(\text{\em $T_Y$ rejects $H_0$}) \leq \alpha. \]
\item[{\em (b)}] Define the following data-dependent subset of $\Theta$:
\begin{equation}
\label{eq:plausibility.region}
C_\alpha(y) = \bigl\{ \vartheta \in \Theta: p_y(\{\vartheta\}) > \alpha \bigr\}.
\end{equation}
Then the frequentist coverage probability of $C_\alpha(Y)$ is lower-bounded by $1-\alpha$, making \eqref{eq:plausibility.region} a nominal $100(1-\alpha)$\% confidence region.  That is, 
\[ \inf_\theta \prob_{Y|\theta}\{C_\alpha(Y) \ni \theta\} \geq 1-\alpha. \]
\end{enumerate}
\end{thm}

This makes clear my previous claim that, if one has a valid inferential model, then decision procedures with good properties are immediate.  It is worth pointing out again that Theorem~\ref{thm:freq}(a) hold for all hypotheses $A$, even those problematic ones concerning a non-linear feature of the full parameter $\theta$ like the ones that cause the false confidence phenomenon in additive inferential models.  Moreover, if $\phi=\phi(\theta)$ is a some feature of interest, then the set 
\[ \bigl\{\varphi \in \phi(\Theta): \textstyle \sup_{\vartheta: \phi(\vartheta) = \varphi} p_y(\{\vartheta\}) > \alpha \bigr\}, \]
is nominal $100(1-\alpha)$\% confidence region for $\phi$.  How to specifically construct a valid inferential model will be discussed in Section~\ref{S:valid.im}.

\section{Remarks}
\label{S:remarks1}

\subsection{BFFs and philosophy of science}
\label{SS:bff}

\citet{xl.bff.2017} describes the efforts of a group working in and around the foundations of statistics and probability.  That group goes by the name of ``BFF'' which has a double-meaning: the first, ``Bayes, Frequentist, and Fiducial,'' is one that only a statistician might guess, while the second, ``Best Friends Forever,'' is a modern pop culture term of endearment.  The first two letters---B and F---are familiar to statisticians and below I briefly summarize these and make a connection to perspectives in philosophy of science.
\begin{itemize}
\item {\em Bayesian}. The defining feature of the Bayesian approach is that it is based strictly on probability calculus.  That is, one starts with a statistical model that relates observable data to some unknown parameters, adds a prior probability distribution that is meant to describe uncertainty about those parameters, and then this prior is updated, via Bayes's formula, to a posterior distribution about the parameter given the observed data.  The rationale behind the Bayesian approach to learning is that, roughly, with sufficiently informative data, perhaps after several iterations of this prior-to-posterior updating, those hypotheses which are ``true'' will have posterior probability close to 1, everything else having posterior probability close to 0.  That is, learning is achieved by confirmation, in the spirit of \citet{carnap1962}, \citet{jeffreys1961}, and others, through the aggregation of evidence supporting a hypothesis.     
\vspace{-2mm}
\item {\em Frequentist}. It is possible to talk about a Bayesian {\em approach} because the framework is normative, i.e., it gives a recipe for carrying out an analysis.  It is much harder, however, to describe a frequentist {\em approach}, it is more of a philosophy or a sense in which any statistical approach or method can be evaluated.  The idea is that, for a relevant scientific hypothesis, statistical methods should be good at assessing the plausibility of that hypothesis, where ``good'' is in the sense that the probability they detect departures between data and the hypothesis is low when there is none and high when there is.  But since not detecting departures between data and a hypothesis does not imply truthfulness of the hypothesis, this strategy cannot be used for direct confirmation, only for refutation.  That is, with sufficiently informative data and sufficiently good statistical methods, any hypothesis that is false is likely to be refuted.  This lines up with the views of \citet{popper1959}, \citet{lakatos1978},  and \citet{mayo.book.1996} who argue that knowledge is gained through falsification.  
\end{itemize}

I have chosen to split the discussion of Bayesian and frequentist into two separate bullet points, but they are not mutually exclusive.  Indeed, since frequentism is not an {\em approach}, but rather a means of evaluating a method, one can study the frequentist properties of Bayesian methods.  It is often the case that Bayesian solutions agree with classical frequentist solutions asymptotically, e.g., a Bernstein--von Mises theorem holds, hence the Bayesian solutions also have the desirable frequentist properties asymptotically.  In cases like these where Bayesian posterior probabilities have certain calibration properties, a verification--falsification balance can be achieved; see \citet{box1980} and \citet{rubin1984}.  However, these niceties are not universal since, as the normal coefficient of variation example, the false confidence theorem, and similar results in \citet{fraser2011} and \citet{fraser.etal.2016} indicate, these asymptotic guarantees may not be fully satisfactory and, as \citet{gleser.hwang.1987} remark, might even be misleading in certain cases.  


Compared to the B and F above, it is more difficult to pin down precisely what the second F, {\em fiducial}, is.  This is partly because, as \citet{zabell1992} elucidates, Fisher's explanation of the ``fiducial argument'' changed significantly over time.  At its genesis, Fisher's idea was to derive a sort of posterior distribution for the unknown parameter without a prior and without invoking Bayes's theorem, an ambitious goal.  To see how this works, let $\theta$ be an unknown scalar parameter and suppose that $T=T(Y)$ is a scalar sufficient statistic, having distribution function $F_\theta$. Take any $p \in [0,1]$ and assume that the equation $p=F_\theta(t)$ can be solved uniquely for $t$, given $\theta$, and for $\theta$, given $t$.  That is, assume there exists $t_p(\theta)$ and $\theta_p(t)$ such that 
\[ p = F_\theta(t_p(\theta)) = F_{\theta_p(t)}(t), \quad \forall \; (t,\theta). \]
If the sampling model is ``monotone'' in the sense that, for all $(p,t,\theta)$, 
\[ t_p(\theta) \geq t \iff \theta_p(t) \leq \theta, \]
then it is an easy calculation to show that 
\begin{equation}
\label{eq:fiducial}
p = \prob_{T|\theta}\{T \leq t_p(\theta)\} = \prob_{T|\theta}\{\theta_p(T) \leq \theta\}. 
\end{equation}
The above expression is a mathematical fact, but has no direct connection to the inference problem because it doesn't depend on the observed $T=t$.  Fisher's big idea was to fix $T=t$ and {\em define} the fiducial probability of the event $\{\theta \geq \theta_p(t)\}$ to be equal to $p$.  In other words, for $t$ the observed value of $T$, Fisher {\em defined} 
\[ \text{``$\prob\{\theta \geq \theta_p(t)\}$''} = p, \quad \forall \; p \in [0,1]. \]
The collection $\{\theta_p(t): p \in [0,1]\}$, for fixed $t$, defines the quantiles of a distribution and, therefore, a distribution itself.  Fisher called this the {\em fiducial distribution} of $\theta$, given $T=t$.  Naturally, since the fiducial probabilities are ``derived'' from calculations based on the sampling distribution, $\prob_{T|\theta}$, of $T$ for fixed $\theta$, Fisher's interpretation of the fiducial probabilities was in a frequentist sense.  So it's no surprise that the distinction between Fisher's fiducial quantile $\theta_p(t)$ and Neyman's $100p$\% upper confidence limit for $\theta$ is blurry.  Fisher ran into trouble, however, in the extension of his fiducial argument to multi-dimensional parameters.  The challenges he faced were, in particular, non-uniqueness and a failure of the fiducial probability for some hypotheses to maintain the calibration property suggested by the connection to the sampling distribution in the expression \eqref{eq:fiducial}.  Notice that these challenges Fisher faced are both part of the high-level message coming out of the false confidence theorem in Section~\ref{SS:fct}, i.e., that data is not informative enough to uniquely identify precise and reliable probabilities for every hypothesis.  Fisher's response to these challenges was to back away from the sampling interpretation of his fiducial probabilities, opting instead for a more subjective, conditional view.  This effort, however, was unsuccessful, which led to fiducial being labeled Fisher's ``one great failure'' \citep[][p.~369]{zabell1992} and ``biggest blunder'' \citep[][p.~105]{efron1998}.  \citet[][p.~382]{zabell1992} summarizes Fisher's motivation and efforts beautifully:
\begin{quote}
{\em Fisher's attempt to steer a path between the Scylla of unconditional behaviorist methods which disavow any attempt at ``inference'' and the Charybdis of subjectivism in science was founded on important concerns, and his personal failure to arrive at a satisfactory solution to the problem means only that the problem remains unsolved, not that it does not exist}.
\end{quote}
There are many other nice papers that discuss Fisher's views and insights on fiducial and inverse probability, among other things; see, e.g., \citet{aldrich2000, aldrich2008}, \citet{edwards1997}, \citet{savage1976}, \citet{seidenfeld1992}, and \citet{zabell1989}.  

Both the aforementioned B and F correspond to well-established schools of thought in the philosophy of science, so what about fiducial?  If the fiducial argument is difficult to understand even as a statistical approach or method, then understanding what, if any, philosophical position it corresponds to is all the more difficult.  Maybe there's a new branch of philosophy of science waiting to be discovered.  A fiducial-style philosophy of science would be similar to a confirmationist's view in the sense that it would be based on probabilities (or some other degree of belief measure) but also similar to a falsificationist's view in the sense that these probabilities can only be used as a tool for refutation, not for proof.  In order to make such a theory work, I think some form of non-additivity would be necessary to handle the important and very common case where data {\em isn't inconsistent} with either $A$ or $A^c$.  I'll leave it to the experts in philosophy of science to decide if the precise formulation of such a theory is worth pursuing.

\subsection{``The most important unresolved problem''}

Even if it's not what Fisher had intended, at a high level, the fiducial argument can be viewed as an attempt to achieve some sort of middle-ground between the more widely-accepted B and F, one that would provide prior-free, data-dependent degrees of belief, like a Bayesian posterior, that simultaneously have certain frequentist calibration properties.  \citet{savage1961} described the fiducial argument as ``a bold attempt to make the Bayesian omelet without breaking the Bayesian eggs.''  As the quote from Zabell above communicates, the fact that Fisher's solution was unsuccessful doesn't mean that a satisfactory solution isn't possible or that the problem itself isn't interesting and worth pursuing.  Indeed, \citet[][p.~41]{efron.cd.discuss} says that 
\begin{quote}
{\em ...perhaps the most important unresolved problem in statistical inference is the use of Bayes theorem in the absence of prior information.} 
\end{quote}
Similar comments can be found in \citet{efron1998} and, even more recently, during his featured address during the 2016 BFF workshop at Rutgers,\footnote{\url{https://statistics.rutgers.edu/bff2016}} Efron referred to the above problem as the {\em Holy Grail} of statistics.  

To be clear, there is no problem at all in getting {\em a} posterior distribution without using a prior or Bayes's theorem.  Any data-dependent probability distribution on the parameter space would be a candidate.  What Efron is talking about is more than just getting a posterior distribution without a prior; he's asking for a ``prior-free'' procedure to construct a ``posterior distribution'' that is reliable, that has performance guarantees.  It's in this sense that my efforts here, which do not have anything directly to do with Bayes's theorem, are relevant to Efron's ``most important unresolved problem.''  

I claim that the key to solving this most important problem is actually very simple: one just needs to be willing to abandon the requirement that degrees of belief be additive.  As the false confidence theorem implies, however you construct your additive degrees of belief, there will exist hypotheses for which the probabilities lack the desired calibration and, consequently, would be ``unacceptable'' according to Reid and Cox.  So, {\em going non-additive}, i.e., adopting a framework in which degrees of belief are measured by some non-additive set function, is the only way to proceed towards a solution to this open problem.  But non-additivity alone isn't enough.  To avoid false confidence---and solve the fundamental problem---an additional constraint is needed, something like the validity condition in Section~\ref{S:valid}.  There is, of course, no free lunch, so any new framework that can achieve this goal must necessarily make sacrifices elsewhere; see Section~\ref{S:discuss}.

\subsection{Beyond $p$ and $P$}
\label{SS:beyond}

In Summer 2018, the statistics community was buzzing about p-values, their interpretation, and their limitations.  Naturally, this got me thinking about the statistical ``debates'' pitting p-values and the NHST (null hypothesis significance testing) camp against posterior probabilities, Bayes factors, and the Bayesian camp.  To set the scene, while p-values have received---and withstood---more than their fair share of scrutiny over the years \citep[e.g.,][]{bergersellke1987, schervish1996}, the discussion heated up when scientists started pointing fingers at the statistical community for low replication rates, etc.  A bold move was made in 2015 when the editors of {\em Basic and Applied Social Psychology}, or {\em BASP}, banned the use of p-values---and much of the modern statistical toolbox---in their journal \citep{pvalue.ban}.  Naturally, this attracted media attention \citep[e.g.,][]{woolston2015} and even sparked an unprecedented response from the American Statistical Association \citep{wasserstein.lazar.asa}, with a follow-up conference\footnote{\url{https://ww2.amstat.org/meetings/ssi/2017/}} and special issue of {\em The American Statistician}.\footnote{\url{https://amstat.tandfonline.com/toc/utas20/73/sup1}} 

My suspicion, which was confirmed by Professor David Trafimow in a personal communication, was that banning p-values was not a solution in and of itself but, rather, an impetus for change, forcing researchers to think harder about how they reason from data to conclusion, potentially leading to new ideas and higher quality research.  Whether this strategy will be successful remains to be seen \citep[e.g.,][]{fricker.etal.ban}, but almost anything is better than researchers ``$p$-ing'' everywhere \citep{valentine.basp}.  Still, it makes sense to ask if contributions from the statistics community have been helpful.  The dust is perhaps still settling, but, based on the expert commentary supplementing the ASA statement, the follow-up by \citet{ionides.response}, the recent proposal of \citet{rss2017} to simply change the default ``$p < 0.05$'' cutoff to ``$p < 0.005$,'' and the plethora of criticism that followed \citep[e.g.,][]{crane.phack, trafimow.etal.rss, mcshane.etal.rss, lakens.etal.rss}, it seems that statisticians' views are as diverse as ever.  

The response from the Bayesian or, at least, the non-NHST camp is that it's better in some sense to rely on probabilities for making inferences.  As part of the justification for the p-value ban, \citet{pvalue.ban} say
\begin{quote}
{\em the [$p$-value] fails to provide the probability of the null hypothesis, which is needed to provide a strong case for rejecting it.} 
\end{quote}
That a p-value is {\em not} the probability that the null hypothesis is true is, of course, a mathematical fact.  But is it really necessary to have such a probability for making inferences?  Fisher says no:
\begin{quote}
{\em As Bayes perceived, the concept of mathematical probability affords a means, in some cases, of expressing inferences from observational data, involving a degree of uncertainty, and of expressing them rigorously...
yet it is by no means axiomatic that the appropriate inferences...
should always be rigorously expressible in terms of this same concept.} \citep[][p.~40]{fisher1973}
\end{quote}
And about the p-value---or observed significance level---he says 
\begin{quote}
{\em It is more primitive, or elemental than, and does not justify, any exact probability statement about the proposition.} \citep[][p.~46]{fisher1973}
\end{quote}
But before asking about whether probability is a necessity, an even more basic question is if a genuine probability is available in general.  If we put together the famously provocative quote from \citet{definetti.vol1}---{\em probability does not exist}---with the refreshingly blunt remark of Fraser's on page~\pageref{quote:fraser.real} above, then we start to see how concerning this is.\footnote{To be fair, \citet{pvalue.ban} acknowledge that the use of posterior distributions based on default priors are not satisfactory substitutes for genuine prior probabilities.}  If probabilities exist only in our minds, and don't become ``real'' simply by updating via Bayes's formula, then how can they be justified as relevant to a scientific problem?  For example, reconsider the highly-publicized predictions of 2016 U.S.~presidential election by Nate Silver on his {\em FiveThirtyEight} website/blog \citep[e.g.,][]{crane.martin.election}.  On the day of the election, Silver predicted Hilary Clinton would win with probability 0.72 and, therefore, Donald Trump would win with probability 0.28.  These numbers are based on specified prior distributions, statistical models, and polling data, all combined via Bayes's formula and evaluated on a computer running Markov chain Monte Carlo.  Does Silver making such a probability statement suddenly imply that the election will be determined by the flip of a 72/28 coin?  Of course not.  The numbers may be meaningful to Silver himself in many ways, e.g., prices for bets he's willing to make, but the relevance of Silver's numbers to me, to any other consumer, or to ``reality'' must be established.  As I see it, there are only two ways this relevance argument can be made: either I understand and agree with all of Silver's modeling assumptions, and accept his subjective probabilities as my own, or I simply trust that Silver's predictions are accurate or otherwise reliable.  The former is quite unlikely, especially for consumers without statistical training; the latter is a choice that I'm free to make or not, though recent results suggest that Silver's predictions may not be so reliable \citep{crane.market, crane.fpp}.  In any case, (posterior) probabilities are not objective or otherwise obviously relevant to the scientific problem at hand, their relevance has to be established in some way.  

While I am skeptical about being able to justify the relevance of a probability distribution to a scientific problem,\footnote{Consider an ideal case where a {\em true} and, therefore, relevant prior exists, e.g., with exchangeable sequences where de~Finetti's theorem applies \citep{definetti1937, hewitt.savage.1955}.  Even in such cases, the true prior, whose existence is guaranteed by the theorem, depends on the data, so one actually needs an infinite---or at least long---data sequence to identify the true prior.} it is worth pointing out that there would be a tremendous payoff if, in fact, it could be done.  Arguably, the more natural way to think about knowledge is as a process by which evidence piles up in favor of some hypothesis until which its truthfulness is undisputed.  This confirmationinst style of learning is quite powerful because one effectively has access to ``truth'' via probability $\approx 0$ or $\approx 1$ events and what \citet{shafer2016.turns40} calls {\em Cournot's Principle}; see also \citet{shafer2007} and \citet{shafer.vovk.2006}.  Popper's falsificationist view emerged, not because it's more powerful, but because, like me, he doubted the confirmationist's premise that a meaningful probability distribution---one that could be continuously updated via Bayes's theorem, eventually pointing to ``truth''---generally exists.  Popper's only option, therefore, is to discard the strong existence assumption, but the price he pays is not being able to judge a hypothesis as ``true,'' only as ``not yet proven false.''  


\ifthenelse{1=1}{}{
\subsection{Skin in the game}
\label{SS:skin}



{\color{red}
We are in a time where statisticians, at least those of us in academia, are not generally the ones engaged in real data analysis.  The primary focus of academic statisticians is the development of {\em statistical methods}, presumably to be used by practitioners and scientists who do analyze real data.  Of course, specialization has certain efficiency advantages, like the assembly line, but there are consequences to this level of detachment.  

The academic statistician's routine is to develop a method, write a paper that highlights the method's practical and theoretical performance compared to its competitors, and create an R software package \citep{Rmanual} that makes it easy for others to apply the method in their own application.  Those who follow this routine benefit in a number of ways: first, there's the credit associated with publishing a paper; second, by making data software available, there's the credit associated with being an active participant in the {\em reproducible research movement}; and third, there's credit associated with citations to the paper and software package by those who develop the next method and need to carry out their comparisons.  What's interesting is that none of these incentives have to do 

I'm not suggesting that the only valuable contribution a statistician can make is one that is

Of course, there's lots more to an effective statistical analysis, and having skin in the game, than valid inferences, but this is one aspect that statisticians, even those not directly involved in the analysis, can help control.  And enforcing a strong validity constraint like that in Section~\ref{S:valid}

The problem is that statistical methods are perceived to have scientific value even if it's never demonstrated.  

It's fine for statisticians to work on things that may not have an immediate scientific impact


Most statisticians will ignore such concerns and continue with business as usual, developing and publishing statistical methods and software implementations, known to have potentially serious limitations, without warning.  This puts users of these methods at risk and the ``use at your own risk'' clause from the developers of these methods creates a harmful asymmetry.  To acknowledge false confidence means putting {\em skin in the game} \citep{taleb.skin} and I'll discuss this more in Section~\ref{SS:skin}.  


...


}
}

\subsection{A Bayesian version of validity?}
\label{SS:bayes.valid}

According to the false confidence theorem, validity is not compatible with a Bayesian inferential model or, for that matter, any other the produces additive degrees of belief as output.  But the validity condition, in Equation \eqref{eq:valid}, is rather strong, so one could ask if the Bayesian approach satisfies perhaps some weaker form of the validity condition.  The answer to this question depends on what kind of Bayesian one is.  

As I see it, there are two kinds of Bayesians: those who believe in the prior and those who don't.  By ``believe in the prior'' I mean that one would actually bet based on their prior probabilities, would effectively dismiss hypotheses with low prior probability, and/or would be satisfied in evaluating the performance of their posterior distribution with respect to the corresponding marginal distribution of $Y$, i.e., $\prob_Y = \int \prob_{Y|\theta} \, \Pi(d\theta)$, where $\Pi$ is the prior distribution for $\theta$.  If one believes in the prior in this sense, then a {\em Bayesian validity property} comes for free.  That is, by Markov's inequality and the martingale property of the posterior with respect to the $\prob_Y$ marginal, 
\begin{equation}
\label{eq:bayes.valid}
\prob_Y\{\Pi_Y(A) > \alpha^{-1} \Pi(A)\} \leq \alpha, \quad \text{for any $\alpha \in (0,1)$ and any $A$}. 
\end{equation}
In words, the posterior probability being a large multiple of the prior probability is a rare event with respect to $\prob_Y$.  Is there anything more than a superficial similarity between the inequalities in \eqref{eq:bayes.valid} and that in \eqref{eq:valid}?  While there's no notion of a hypothesis about $\theta$ being ``true'' or ``false'' in this context, the closest thing is to consider those hypotheses with small prior probability, since a believer in their prior probabilities would consider these to be ``effectively false.''  Then \eqref{eq:bayes.valid} says that it's rare to see data that provide strong support for hypotheses having small prior probability, which is a version of the validity property.  These self-consistency properties are expected, that is, if the prior is meaningful, then the posterior obtained by the usual probability calculus is equally meaningful.  

The second group---those who don't believe in the prior---consists of those Bayesians who might, for simplicity, use default or non-informative priors, trusting that the prior-to-posterior updates will correct for any shortcomings in the prior specification.  In order for the marginal distribution and the prior probabilities in \eqref{eq:bayes.valid} to make sense, let me assume that the prior $\Pi$ is proper but diffuse, as would often be the case in situations where genuine prior information is lacking.  It's possible that $\prob_Y$ could be very different from the true distribution of $Y$, in which case \eqref{eq:bayes.valid} isn't meaningful at all.  But let me also assume that $\prob_Y$ is a reasonable approximation of the true distribution, e.g., it passes certain prior- or posterior-predictive checks \citep[e.g.,][]{gelman.meng.stern.1996, box1980, rubin1984}.  Even in this case, there is still trouble coming from \eqref{eq:bayes.valid}.  The inequality says that the posterior probabilities are unlikely to be much larger than the prior probabilities.  But virtually any ``non-extreme'' hypothesis $A$ has small $\Pi(A)$ by virtue of the diffuse non-informative prior, {\em not because there's reason to believe that it's false}.  So if a potentially true hypothesis has small prior probability, and it's unlikely that the posterior will be much larger than the prior, then how can one expect to learn?  This is precisely what's going on in the satellite conjunction analysis example summarized in Section~\ref{SS:conjunction}.  The diffuse prior behind the scenes assigns small prior probability to the collision event and, even if collision is imminent, \eqref{eq:bayes.valid} implies that the posterior probabilities are unlikely to turn large, hence false confidence.  So, the inequality \eqref{eq:bayes.valid} that is an asset to the Bayesian who believes in his/her prior becomes a liability to the Bayesian who doesn't.


\section{Non-additive beliefs}
\label{S:nonadditive}

\subsection{Lower probabilities and capacities}

In the previous sections I have argued that, in order to avoid false confidence, it is necessary that the inferential model output be {\em non-additive}, e.g., that $b_y(A) + b_y(A^c) \neq 1$ for some hypothesis $A$.  This section describe some particular non-additive set functions that have been suggested in the literature as measures of degrees of belief, all of which generalize the familiar additive probability.  (Here it is not necessary to consider beliefs that depend on data $y$ so I will write $b$ instead of $b_y$ throughout.)  

The most general non-additive measure considered here is a {\em lower probability}.  Let $\Theta$ be completely metrizable and consider a collection $\M$ of probability distributions on the measurable space $(\Theta, \B_\Theta)$.  Now define the lower probability as 
\[ b(A) = \inf_{P \in \M} P(A), \quad A \in \B_\Theta. \]
Of course, there is a dual to $b$, namely, 
\[ p(A) = \sup_{P \in \M} P(A) = 1-b(A^c), \]
which in the this context is called the corresponding {\em upper probability}.  The lower probability is non-additive in the sense above or, more specifically, they are {\em super-additive} in the sense that  
\[ b(A \cup B) \geq b(A) + b(B) - b(A \cap B), \]
with strict inequality at least for some $A$ and $B$.  Super-additivity implies that $b(A) \leq p(A)$, and this gap between the lower and upper probabilities can be understood as a measure of how imprecise the imprecise model is; \citet{dempster2008} referred to the difference $p(A)-b(A)$ as the ``don't know'' probability.  \citet{walley1991} presents a theory of imprecise probability based on lower probabilities---or, more generally, {\em lower provisions}---with a focus on achieving a {\em coherence}, i.e., that a gambler who determines the prices she's willing to pay to make certain bets according to $b$ can't be made a sure loser.  In fact, if $\M$ is closed and convex, then the lower probability $b$ is coherent \citep[e.g.,][]{walley2000}.  

A special case of lower probabilities that often arises in statistical applications, namely, in robustness studies, are {\em Choquet capacities} or, {\em capacities} for short \citep{choquet1953}.  Let $b$ be a set function defined on a domain $\A \subseteq 2^\Theta$, with the basic properties $b(\varnothing) = 0$ and $b(\Theta)$.  Then $b$ is a {\em capacity of order $n$} if 
\[ b(A) \geq \sum_{\varnothing \neq I \subseteq \{1,\ldots,n\}} (-1)^{|I|-1} b\Bigl( \bigcap_{i \in I} A_i \Bigr) \]
for every $A \in \A$ and every collection $\{A_1,\ldots,A_n\}$ of sets in $\A$ with $A_i \subset A$.  Such a function $b$ is said to be {\em $n$-monotone}.  Then the corresponding dual, $p(A) = 1-b(A^c)$, is {\em $n$-alternating} and satisfies 
\[ p(A) \leq \sum_{\varnothing \neq I \subseteq \{1,\ldots,n\}} (-1)^{|I|-1} p\Bigl( \bigcup_{i \in I} A_i \Bigr). \]
Since capacities are special cases of lower probabilities/provisions, it follows that they too are super-additive.  It is also clear that every capacity of order $n$ is also a capacity of order $n-1$ so the class is shrinking in $n$.  Therefore, capacities of order 2 are the most general and have been investigated by statisticians in, e.g.,  \citet{huber.strassen.1973}, \citet{huber1973.capacity, huber1981}, \citet{wasserman1990}, and \citet{wasserman.kadane.1990}.  

A relevant question is if there are any obvious practical advantages to super-additivity?  The answer is yes, and the simplest is that of modeling {\em ignorance}.  Despite being somewhat extreme, ignorance is apparently a practically relevant case since the developments of default priors are motivated by examples where the data analyst is ignorant about the parameter, i.e., no genuine prior information is available.  So, understanding how to describe ignorance mathematically ought to be insightful for understanding the advantages and limitations of prior-free probabilistic inference.  A first question is: {\em what does ignorance mean?}  A standard assumption is that the universe $\Theta$ is given---hence not fully ignorant---and contains at least two points.  Subject to these conditions, ignorance means that there is no evidence available that supports the truthfulness or falsity of any non-trivial hypothesis, and a natural way to encode this mathematically is by setting $b(A) = 0$ for any proper subset $A$ of $\Theta$, and $b(\Theta)=1$.  \citet{shafer1976} and others call this a {\em vacuous belief}, and it's easy to check that this $b$ is super-additive.  Compare this to standard approaches to formulating ignorance in classical probability, based on assumptions of symmetry/invariance or Laplace's ``principle of indifference,'' both of which are based on certain judgments---or {\em knowledge}---and, therefore, don't really model ignorance.  

When working with ordinary probabilities, updating via conditioning is straightforward and uncontroversial.  But the analogous updating of lower probabilities, capacities, and even the belief functions discussed in Section~\ref{SS:belief} can have some arguably unexpected behaviors, e.g., {\em dilation}.  See \citet{gong.meng.update} for a modern discussion of these updating rules and their respective properties.  Some comments on Dempster's rule of combination \citep[see, e.g.,][]{shafer1976} in the context of combining inferential models represented as belief functions will be given in Section~\ref{SS:efficiency} below.

\subsection{Belief functions}
\label{SS:belief}

A {\em belief function} is a special type of capacity, one that is $\infty$-monotone, i.e., $n$-monotone for all $n \geq 1$.  This is a strong condition, but there are some advantages that come with this specialization, namely, that capacities are quite complicated but belief functions can be characterized in terms of simpler and more intuitive mathematical objects.  

As is customary in the literature on belief functions, I'll assume for the moment that $\Theta$ is a finite space.  Then the power set $2^\Theta$ is also finite and I can imagine a probability mass function $m$ defined on $2^\Theta$, assumed to satisfy $m(\varnothing) = 0$.  For $A \subseteq \Theta$, the quantity $m(A)$ encodes the degree of belief in the truthfulness of $A$ but in no proper subset of $A$.  Then the function 
\begin{equation}
\label{eq:basic.belief}
b(A) = \sum_{B: B \subseteq A} m(B), \quad A \subseteq \Theta,
\end{equation}
is a belief function, i.e., is a capacity of order $\infty$.  Intuitively, the belief in the truthfulness of the assertion $A$ is defined as the totality of the degrees of belief in $A$ and all its proper subsets.  It turns out that there is a one-to-one correspondence between the mass function $m$ and the belief function $b$; see \citet{shafer1976} or \citet{yagerliu2008} for details.  

A further insight that can be gleaned from \eqref{eq:basic.belief} is that the belief function is effectively characterized by a {\em random set}.  That is, one can think of $m$ as the mass function for a random element taking values in $2^\Theta$.  But this equivalence between belief functions and random sets only holds up in the case when $\Theta$ is a finite set.  When $\Theta$ is infinite, it turns out that random sets corresponding to a special class of belief functions, which I describe in the next subsection.  To make the connection between random sets and belief functions in the general case, there is an extra layer, what \citet{shafer1979} calls an {\em allocation of probability}, which I will not discuss here.  \citet{kohlas.monney.hints} also discuss this gap between general belief functions and those based on random sets. 

My discussion of belief functions here does nowhere near justice to the tremendous amount of incredible work that has been done in this area.  And since I can't possibly provide an adequate overview of the literature, let me just point the reader to the 2016 special issue\footnote{\url{https://www.sciencedirect.com/journal/international-journal-of-approximate-reasoning/special-issue/10BG01ZSM7P}} of the {\em International Journal of Approximate Reasoning}, edited by Thierry Den{\oe}ux, to celebrate 40 years since publication of Glenn Shafer's monograph, which gives a look back at the belief function developments over that time.

\subsection{Random sets}
\label{SS:random.sets}

As eluded to in the previous subsection, random sets lead to belief functions but not all belief functions correspond to random sets; the equivalence only holds when the domain $\Theta$ is a finite set.  In general, a random set on an infinite domain will inherit certain continuity properties that a belief function need not satisfy.  

As expected, defining random sets in a mathematically rigorous way requires considerable care.  Since this amount of rigor isn't necessary---or helpful---for my present purposes, I will not be completely precise in what follows.  The reader interested in the precise details can check out \citet{nguyen1978, nguyen.book} or \citet{molchanov2005}.  Let $\T$ denote a random set in $\Theta$, i.e., a (measurable) function from a sample space to $2^\Theta$.  If $\prob_\T$ is the distribution of $\T$, define the {\em containment functional} 
\[ c(A) = \prob_\T(\T \subseteq A), \quad A \subseteq \Theta. \]
This containment functional obviously satisfies $c(\varnothing) = 0$ and $c(\Theta) = 1$.  It can also be shown that $c$ is $\infty$-monotone, so $c$ is a belief function.  However, through its connection to a genuine probability measure, $c$ is also continuous in the sense that 
\[ c(A_n) \downarrow c(A), \quad \text{for any $A_n \downarrow A$}, \]
a condition which is not required of general belief functions.  Interestingly, the additional layer that Shafer needed to connect belief functions to random sets is entirely due to the latter being continuous.  Indeed, according to the famous {\em Choquet capacity theorem}, every continuous belief is the capacity functional of a random set.  In Section~\ref{S:valid.im}, I will describe the construction of a valid inferential model whose output is a random set containment function, hence is equivalent to a framework built around continuous belief functions.  

An even more specific class of belief functions is that based on random sets, $\T$ which are {\em nested} in the sense that, for any pair of sets $T$ and $T'$ in the support of $\T$, either $T \subseteq T'$ or $T \supseteq T'$.  Then the corresponding belief function, $b$, is called {\em consonant}, and the dual, $p$, to $b$ has the following surprising property:
\[ p(A) = \sup_{\vartheta \in A} p(\{\vartheta\}). \]
That is, the plausibility and, hence, the belief function are determined by an ordinary point function, called the {\em plausibility contour}, 
\[ p(\{\vartheta\}) = \prob_\T(\T \ni \vartheta), \quad \vartheta \in \Theta. \]
The belief/plausibility functions corresponding to nested random sets are also known as {\em necessity/possibility measures} \citep[e.g.,][]{dubois.prade.book}.

\section{Valid inferential models}
\label{S:valid.im}

\subsection{Basic construction}

The previous sections have focused primarily on the point that non-additive degrees of belief satisfying a validity constraint (or something similar) are needed in order to avoid false confidence.  Obvious questions that come to mind are {\em does a valid inferential model exist?} and, if so, {\em how to construct it?}  This section reviews the detailed construction first presented in \citet{imbasics} and then expanded upon in \citet{imbook}.  

The inferential model depends explicitly on the posited statistical model, so the starting point for the following construction is a particular representation of that model, as a data-generating process.  That is, write
\begin{equation}
\label{eq:assoc}
Y = a(\theta, U), \quad U \sim \prob_U, 
\end{equation}
where $U \in \UU$ is called an {\em auxiliary variable} with distribution $\prob_U$ fully known, and $a: \Theta \times \UU \to \YY$ is a known mapping.  The representation in \eqref{eq:assoc} is called an {\em association} between data $Y$, parameter $\theta$, and auxiliary variable $U$; this functional form of the statistical model is one of the extra ``$\ldots$'' inputs in the inferential model definition.  An association always exists since any sampling model that can be simulated on a computer must have a representation \eqref{eq:assoc}, but that representation is not unique.  The association need not be based on knowledge of the actual data-generating process, it can be viewed simply as a representation of the data analyst's uncertainty; that is, no claim that $U$ exists in ``real life'' is made or required.  The advantage to introducing $U$ in \eqref{eq:assoc} is that it identifies a fixed probability space, namely, $(\UU, \mathscr{B}, \prob_U)$, for a $\sigma$-algebra $\mathscr{B}$ containing subsets of $\UU$, on which we can carry out probability calculations independent of the data and parameter; as I describe below, these probability calculations are with respect to the distribution of a random set on $\UU$.  

Readers familiar with Fisher's fiducial argument, or the structural inference framework of \citet{fraser1968}, will recognize the auxiliary variable $U$ as a ``pivotal quantity.''  Roughly speaking, these approaches aim to solve the equation \eqref{eq:assoc} for $\theta$ and then use the distribution $U \sim \prob_U$ and the observed value of $Y$ to induce a distribution for $\theta$. This fiducial argument will lead to an additive probability distribution on $\Theta$ which does not satisfy the validity property.  Even the extended fiducial argument put forth by Dempster, where the inferential output need not be additive, the validity property generally fails.  Therefore, something extra is needed.  More discussion on the connection between the fiducial argument and the approach presented here is given in Section~\ref{SS:fiducial}.

According to the formulation in \eqref{eq:assoc}, there is a quantity $U$ that, together with $\theta$, determines the value of $Y$.  If both $Y$ and $U$ were observable, then \eqref{eq:assoc} could be simply solved for $\theta$, hence the best possible inference.  That is, if $Y=y$ and $U=u$ were observed, then it would be {\em certain} that the true $\theta$ belongs to the set 
\[ \Theta_y(u) = \{\vartheta \in \Theta: y = a(\vartheta, u)\}. \]
However, since $U$ is not observable, the best one can do is to accurately \emph{predict}\footnote{Perhaps {\em guess} is a better word to describe this operation than {\em predict} but the former has a particularly non-scientific connotation, which I'd like to avoid.} the unobserved value using some information about its distribution $\prob_U$.  

To accurately predict the specific unobserved value of $U$, it is not enough to use an independent draw from $\prob_U$ as the fiducial argument suggests.  Intuitively, hitting a target with positive probability requires stretching that draw from $\prob_U$ into a genuine set.  Towards formalizing this notion of ``stretching a draw from $\prob_U$,'' consider a random set, $\S$, having distribution $\prob_\S$, with the property that it contains a large $\prob_U$-proportion of $U$ values with high probability.  More precisely, define $\gamma_\S(u) = \prob_\S(\S \ni u)$, a feature of the distribution $\prob_\S$, called the contour function of $\S$, and suppose that 
\begin{equation}
\label{eq:prs}
\text{$\gamma_\S(U)$ is stochastically no smaller than $\unif(0,1)$.}
\end{equation}
In other words, $\S$---or, rather, $\prob_S$---must be such that 
\[ \prob_U\{\gamma_\S(U) \leq \alpha\} \leq \alpha, \quad \forall \; \alpha \in [0,1]. \]
Intuitively, this condition suggests that $\S$ is ``reliable'' in the sense that it can cover a draw $U$ from $\prob_U$ with high probability; more details about the random set and condition \eqref{eq:prs} in Section~\ref{SS:validity.theorem}.  The distribution $\prob_\S$ of the random set $\S$ is the second extra ``$\ldots$'' input in the inferential model definition.  

Given that $\S$ provides a reliable guess of where the unobserved $U$ resides, in the sense of \eqref{eq:prs}, then an equally reliable guess of where the unknown $\theta$ resides, given $Y=y$, is 
\[ \Theta_y(\S) = \bigcup_{u \in \S} \Theta_y(u), \]
another random set, with distribution determined by $\prob_\S$.  It is, therefore, natural to say that data $y$ supports the claim ``$\theta \in A$'' if $\Theta_y(\S)$ is a subset of $A$.  Then the inferential model output is
\begin{equation}
\label{eq:belief}
b_y(A) = \prob_\S\{\Theta_y(\S) \subseteq A\}, \quad A \subseteq \Theta.
\end{equation}
This is a belief function, like in Section~\ref{SS:belief}, so it inherits those properties, e.g., it's super-additive and has a corresponding plausibility function $p_y(A) = 1-b_y(A^c)$.  

It will typically be the case that $\Theta_y(\S)$ is {\em nested}, as discussed in Section~\ref{SS:random.sets}; I'll talk about why in Section~\ref{SS:efficiency}.  Consequently, the belief function $b_y$ is consonant and fully determined by the plausibility contour $p_y(\{\vartheta\}) = \prob_\S\{\Theta_y(\S) \ni \vartheta\}$, $\vartheta \in \Theta$.  According to \citet[][Ch.~10]{shafer1976}, while the use of a consonant belief function as a description of evidence might be questionable in general, one case where it would be appropriate is in describing inferential evidence, as is the focus here.  

Next are a couple of simple examples to illustrate the inferential model construction and how the corresponding belief/plausibility functions are evaluated.  

\begin{figure}[t]
\begin{center}
\subfigure[Plausibility contour]{\scalebox{0.59}{\includegraphics{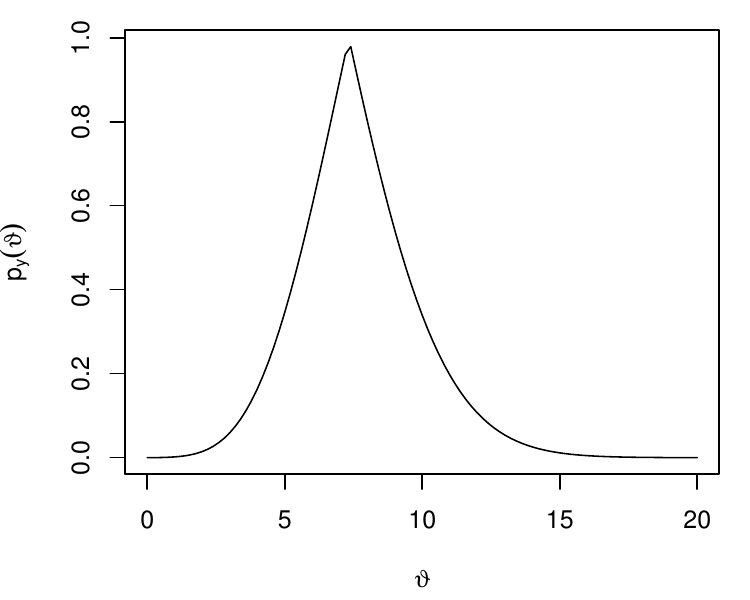}}} 
\subfigure[$b_y(A_\vartheta)$ and $p_y(A_\vartheta)$]{\scalebox{0.59}{\includegraphics{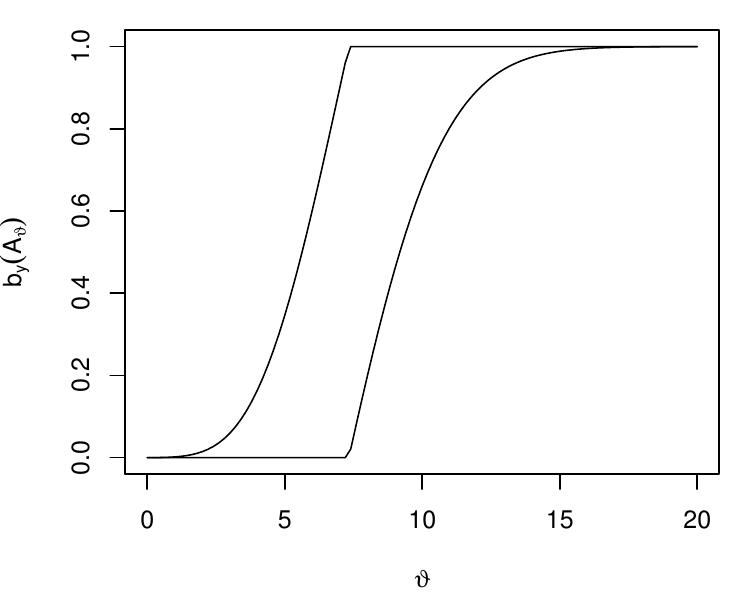}}}
\end{center}
\caption{Results of the gamma inference problem in Example~\ref{ex:gamma} with $Y=7$.  Panel~(a) shows the plausibility contour and Panel~(b) shows the belief and plausibility functions for the collection of hypotheses $A_\vartheta=(0, \vartheta)$.}
\label{fig:gammapl}
\end{figure}

\begin{example}
\label{ex:gamma}
Let $Y$ be a scalar with continuous distribution function $F_\theta$, depending on a scalar parameter $\theta$.  Then a natural choice for the association \eqref{eq:assoc} would be $Y = F_\theta^{-1}(U)$, where $U \sim \unif(0,1)$.  Then $\Theta_y(u) = \{\vartheta: F_\vartheta(y) = u\}$ and, with a random set $\S$, it is easy to check that $\Theta_y(\S) = \{\vartheta: \S \ni F_\vartheta(y)\}$, so 
\[ p_y(\{\vartheta\}) = \prob_\S\{\S \ni F_\vartheta(y)\}. \]
For instances like this, where the auxiliary variable $U$ is $\unif(0,1)$, a good default choice of the random set is 
\begin{equation}
\label{eq:default.prs}
\S = \{u \in (0,1): |u-0.5| \leq |U-0.5|\}, \quad U \sim \unif(0,1). 
\end{equation}
This is just a sub-interval of $(0,1)$, centered at 0.5, with a random width.  For this choice of $\S$, the plausibility contour above can be simplified:
\[ p_y(\{\vartheta\}) = 1 - |2 F_\vartheta(y) - 1|. \]
From here, plausibility of any hypothesis $A$ can be evaluated as 
\[ p_y(A) = \sup_{\vartheta \in A} \pl_y(\{\vartheta\}), \]
and hence belief can be evaluated according to the formula $b_y(A) = 1-p_y(A^c)$.  For example, suppose that $F_\theta$ is the distribution function of $\gam(\theta, 1)$, where $\theta > 0$ is the unknown shape parameter.  Panel~(a) in Figure~\ref{fig:gammapl} shows the plausibility contour and Panel~(b) shows $b_y(A_\vartheta)$ and $p_y(A_\vartheta)$ where $A_\vartheta = (0,\vartheta)$, sort of lower and upper distribution functions.  Note the gap between the lower and upper in Panel~(b).  
\end{example}

\begin{figure}[t]
\begin{center}
\subfigure[Plausibility contour]{\scalebox{0.59}{\includegraphics{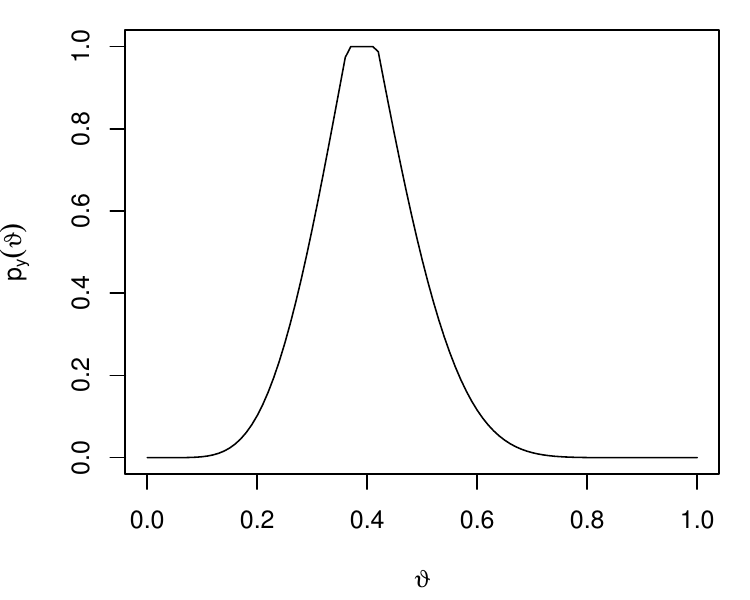}}} 
\subfigure[$b_y(A_\vartheta)$ and $p_y(A_\vartheta)$]{\scalebox{0.59}{\includegraphics{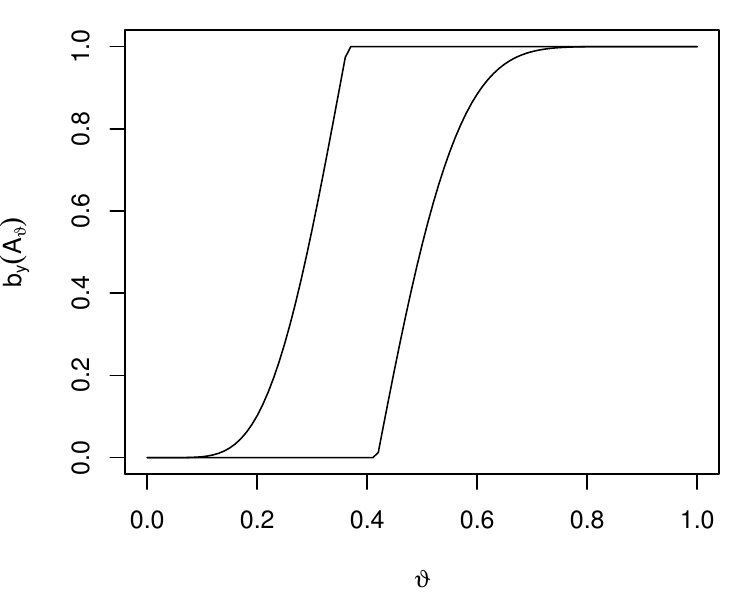}}}
\end{center}
\caption{Results of the binomial inference problem in Example~\ref{ex:binom1} with $y=7$ and $n=18$.  Panel~(a) shows the plausibility contour and Panel~(b) shows the belief and plausibility functions for the collection of hypotheses $A_\vartheta=(0, \vartheta)$.}
\label{fig:binom1}
\end{figure}

\begin{example}
\label{ex:binom1}
Consider a binomial model, $Y \sim \bin(n,\theta)$, where $n$ is a known positive integer and $\theta \in (0,1)$ is the unknown success probability.  Despite its apparent simplicity, inference on $\theta$ in the binomial model is far from trivial \citep[e.g.,][]{bcd2001}.  In this case, the most natural version of \eqref{eq:assoc} is 
\[ F_\theta(Y-1) \leq U < F_\theta(Y), \quad U \sim \unif(0,1), \]
where $F_\theta$ is the $\bin(n,\theta)$ distribution function.  There is no simple equation linking $(Y, \theta, U)$ in this case, just a rule ``$Y=a(\theta,U)$'' for producing $Y$ with given $\theta$ and $U$.   This determines the map $u \mapsto \Theta_y(u)$, which is given by 
\begin{align*}
\Theta_y(u) & = \{\vartheta: F_\vartheta(y-1) \leq u < F_\vartheta(y)\} \\
& = \{\vartheta: G_{n-y+1,y}(1-\vartheta) \leq u < G_{n-y,y+1}(1-\vartheta)\} \\
& = \{\vartheta: 1-G_{n-y+1,y}^{-1}(u) \leq \vartheta < 1-G_{n-y,y+1}^{-1}(u)\}, 
\end{align*}
where $G_{a,b}$ denotes the ${\sf Beta}(a,b)$ distribution function.  This calculation follows from the connection between the binomial and beta distribution functions, namely, $F_\theta(y) = G_{n-y, y+1}(1-\theta)$, which can be derived using integration-by-parts.  Note that, in this discrete-data setting, the set $\Theta_y(u)$ is an interval, not a singleton like in the continuous-data example above.  Next, for simplicity, again I'll use the default random set $\S$ in \eqref{eq:default.prs} for predicting the unobserved $U \sim \unif(0,1)$. Combining $\S$ with $\Theta_y(\cdot)$ above gives 
\[ \Theta_y(\S) = \bigcup_{u \in \S} \Theta_y(u) = \bigl[ 1 - G_{n-y+1, y}^{-1}(\tfrac12 + |\tilde U - \tfrac12|), 1 - G_{n-y, y+1}^{-1}(\tfrac12 - |\tilde U - \tfrac12|) \bigr], \]
where $\tilde U \sim \unif(0,1)$.  Then the probability, $\prob_\S\{\Theta_y(\S) \subset A\}$ can be computed for any $A$ using the distribution of the endpoints of $\Theta_y(\S)$, which are simple functions of a uniform random variable.  For example, the plausibility contour equals 
\[ p_y(\{\vartheta\}) = \prob_\S\{\Theta_y(\S) \ni \vartheta\} = 1 - \prob_\S\{\max \Theta_y(\S) < \vartheta\} - \prob_\S\{\min \Theta_y(\S) > \vartheta\}. \]
Figure~\ref{fig:binom1}(a) shows a plot of this plausibility contour based on and observed $y=7$ and $n=18$.  A feature that distinguishes this plausibility contour from that in Example~\ref{ex:gamma} is the plateau at the top; this is a consequence of $\Theta_y(u)$ being an interval for all $u$, which itself is a consequence of the discreteness of $Y$.  Panel~(b) shows the belief and plausibility for a sequence of one-sided hypotheses like that in Figure~\ref{fig:gammapl}(b).  
\end{example}

\begin{figure}[t]
\begin{center}
\subfigure[$\|y\|^2 = 0.5$ vs $t=1$]{\scalebox{0.6}{\includegraphics{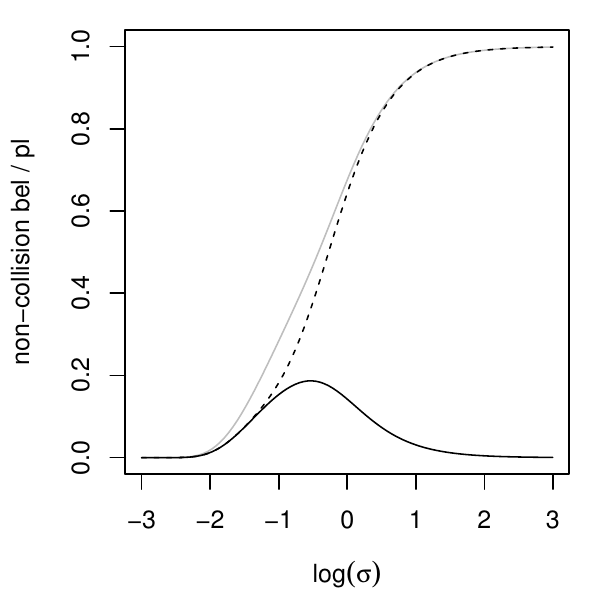}}}
\subfigure[$\|y\|^2 = 0.9$ vs $t=1$]{\scalebox{0.6}{\includegraphics{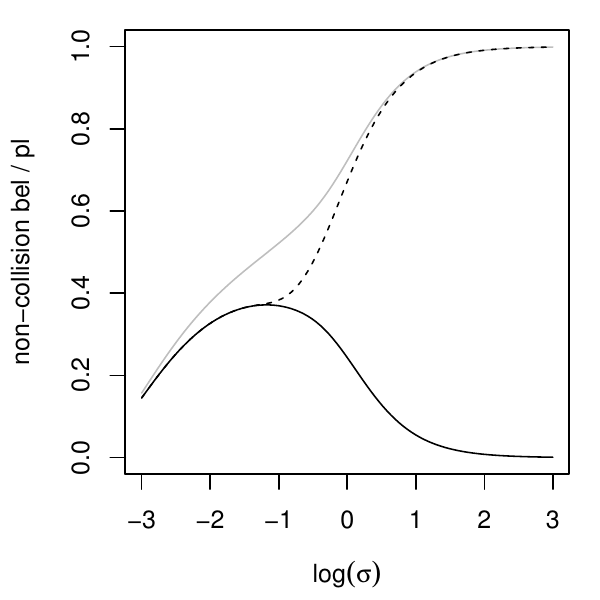}}}
\subfigure[$\|y\|^2 = 1.1$ vs $t=1$]{\scalebox{0.6}{\includegraphics{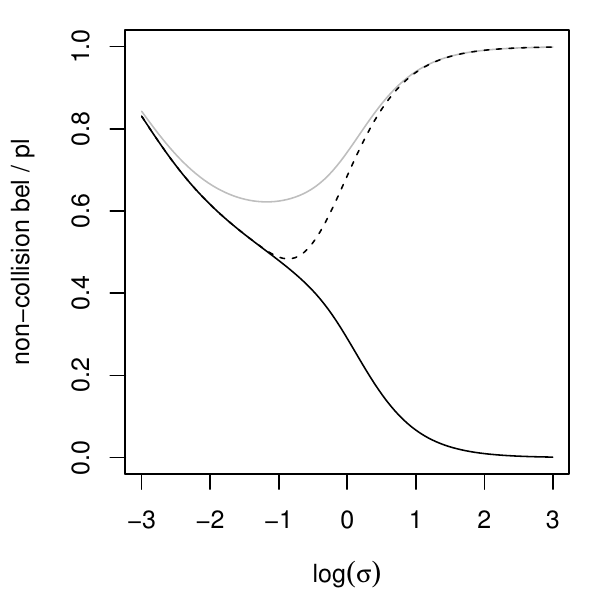}}}
\subfigure[$\|y\|^2 = 1.5$ vs $t=1$]{\scalebox{0.6}{\includegraphics{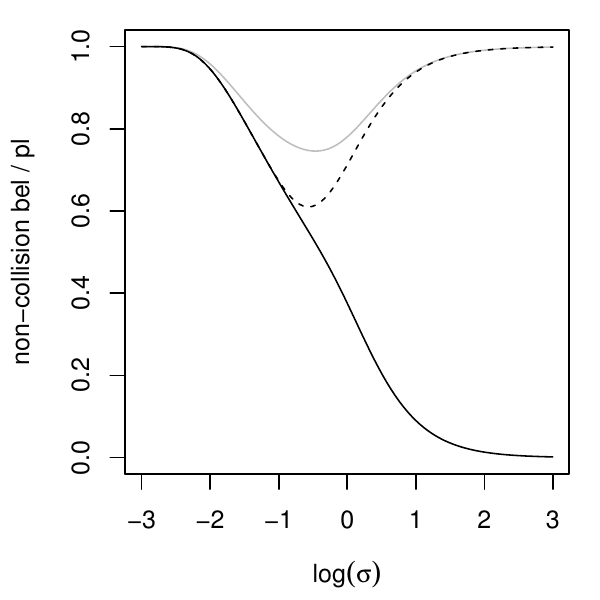}}}
\end{center}
\caption{Like in Figure~\ref{fig:collision}, plots of the non-collision belief (solid) and plausibility (dashed) based on a valid inferential model as a function of $\sigma$ for several configurations of $y$.  For reference, the non-collision probability from \eqref{eq:ncp} is shown in gray.}
\label{fig:collision2}
\end{figure}

To conclude this subsection, I briefly revisit the satellite conjunction analysis example from Section~\ref{SS:conjunction}.  Details about the calculations will be presented elsewhere, here I only want to focus on the non-additivity feature and its interpretation.  Figure~\ref{fig:collision2} considers exactly the situation investigated in Figure~\ref{fig:collision}, but displays the belief and plausibility in the non-collision hypothesis, based on a valid---and, in fact, ``optimal''---inferential model, as a function of the error variance $\sigma$, for several different estimates, $y$, of the distance between the two satellites; for reference, the Bayesian non-collision probabilities from \eqref{eq:ncp} are displayed in gray.  First note that, in the ideal, high-quality data case with small $\sigma$, the belief and plausibility values are equal and generally close to the Bayesian non-collision probability.  But as $\sigma$ increases and the quality of the data is diminished, the non-collision probability is forced up to 1 but the non-collision belief and plausibility values split and eventually reach 0 and 1, respectively.  What does this mean?  \citet{dempster2008} calls the gap between belief and plausibility a ``don't know probability,'' and that's the appropriate way to interpret the phenomenon displayed in these plots.  That is, when data is low-quality, a clear determination about the satellites' safety isn't possible, so it's more appropriate to say ``I don't know'' than to conclude that the satellite is safe as the non-collision probability would suggest.  This reveals the advantage of having the extra flexibility coming from non-additivity; that is, having the extra ``don't know'' category implies that data not supporting collision need not support non-collision, which apparently is the key to avoiding false confidence.  

\subsection{Validity theorem}
\label{SS:validity.theorem}

Importantly, the inferential model obtained by the above construction, driven by the random set $\S \sim \prob_\S$ that satisfies \eqref{eq:prs}, is provably valid in the sense of Definition~\ref{def:valid}.  

\begin{thm}[\citealt{imbasics}]
\label{thm:valid}
Suppose that $\S$ satisfies \eqref{eq:prs} and that $\Theta_y(\S)$ is non-empty with $\prob_\S$-probability~1 for each $y$.  Then the inferential model with output $b_y$ defined in \eqref{eq:belief} is valid in the sense of Definition~\ref{def:valid}.
\end{thm}


Recall that the validity condition was designed to properly calibrate the data analyst's beliefs for scientific application, to avoid false confidence.  An important practical consequence, presented in Theorem~\ref{thm:freq} above, is that it's straightforward to construct decision procedures with guaranteed frequentist error rate control.  For example, the corresponding plausibility region \eqref{eq:plausibility.region} attains the nominal frequentist coverage probability.  

Let's discuss the conditions of the theorem, starting with \eqref{eq:prs}.  In particular, how to construct a random set $\S \sim \prob_\S$ such that \eqref{eq:prs} holds?  \citet{imbasics} show that this is actually a relatively mild condition and provide some general strategies for constructing such a set $\S$.  For example, if $\prob_U$ is non-atomic and $h: \UU \to \RR$ is a continuous function, non-constant on any set with positive $\prob_U$-measure, then $h(U)$ is a continuous random variable and the random set 
\begin{equation}
\label{eq:h.prs}
\S := \{u \in \UU: h(u) \leq h(\tilde U)\}, \quad \tilde U \sim \prob_U, 
\end{equation}
with distribution $\prob_\S$ determined by the push-forward of $\prob_U$ through the set-valued mapping above, satisfies \eqref{eq:prs}; see \citet[][p.~61]{imbook} for a proof.  The default random set in \eqref{eq:default.prs} is of this type, with $\UU=[0,1]$, $\prob_U=\unif(0,1)$, and mapping $h(u) = |u-\frac12|$.  Other examples of random sets satisfying \eqref{eq:prs} are presented in Section~\ref{SS:examples}.  

Next, the non-emptiness assumption in the theorem often holds automatically for natural choices of $\S$, but can fail if the dimensions of $\theta$ and $U$ don't match up (see Section~\ref{SS:dimred}) or in problems that involve non-trivial parameter constraints.  When $\prob_\S\{\Theta_x(\S) \neq \varnothing\} < 1$, \citet{leafliu2012} propose to use an {\em elastic} version of $\S$ that stretches just enough that $\Theta_y(\S)$ is non-empty, and they prove a corresponding validity theorem.    

To conclude this brief introduction, I should mention that the approach described above is not necessarily the only way to construct a valid inferential model.  I'm not aware of any competing constructions, but I make no claims that what I've presented above is the best or only way to achieve the validity property.

\subsection{Dimension reduction}
\label{SS:dimred} 

The basic construction presented above is actually quite general, but depending on the data structure and/or the objectives of the analysis, it may be possible to make the inferential model more efficient by reducing the dimension of the auxiliary variable that links the data and parameter in \eqref{eq:assoc}.  This dimension reduction strategy has proved to be useful in two frequently encountered scenarios: 
\begin{itemize}
\item when data consists of individual components each contributing some information about a common parameter, e.g., like with iid data, and
\vspace{-2mm}
\item when the parameter of interest is some lower-dimensional feature, $\phi=\phi(\theta)$, of the full parameter $\theta$.  
\end{itemize}
I'll give a brief summary of the general ideas here, saving details for the examples in Section~\ref{SS:examples}; see, also, \citet[][Ch.~6--7]{imbook} and the relevant references.  

Consider the case of iid data, $Y=(Y_1,\ldots,Y_n)$.  There are two obvious ways to proceed.  The first is to develop an association of the form \eqref{eq:assoc} with $U=(U_1,\ldots,U_n)$, an auxiliary variable for every data point, introduce a random set for the $n$-dimensional $U$, and proceed as above.  The second is to construct valid inferential models with output $b_{y_i}$ for each $i=1,\ldots,n$, and than combine---since they are all independent and carry information about the same $\theta$---combine them in some way, e.g., using Dempster's rule of combination.  It turns out that neither of these two obvious approaches are most efficient; see Section~\ref{SS:efficiency}.  The key insight, first presented in \citet{imcond}, is that if information about $\theta$ is coming from multiple sources, then there will features of the unobservable auxiliary variable $U=(U_1,\ldots,U_n)$ that actually {\em are observed}, hence no reason to predict those features with a random set.  The strategy is to express $U$ in a coordinate system in wherein some coordinates are observed and the others are not.  This is done by identifying a one-to-one mapping $u \mapsto (\tau(u), \eta(u))$, where $\eta(u)$ corresponds to the observable part and $\tau(u)$ the unobservable part.  Then $\tau(U)$ will necessarily be of lower dimension than $U$ itself, so predicting $\tau(U)$ with a suitable random set will be easier and more efficient than predicting the full $U$.  There are a number of strategies available to identify the $(\tau, \eta)$ mapping---using sufficient statistics, symmetry/invariance, a conditional representation of the statistical model, solving a differential equation, etc.---and I will illustrate these different approaches in the examples that follow.  

Next, assuming that the combination steps have been carried out, so that the auxiliary variable $U$ is currently of the same dimension as that of the parameter $\theta$, if the goal is inference on $\phi = \phi(\theta)$, then there may be an opportunity to reduce the dimension of the auxiliary variable further.  In the simplest case, imagine that some feature of $U$ is directly connected to data and $\phi$, while the remaining features of $U$ don't directly connect to $\phi$.  \citet{immarg} show that one can effectively {\em ignore} the part of $U$ unrelated to $\phi$, and predict only that relevant feature of $U$.  Again, that feature would be of dimension lower than that of the original $U$, hence an opportunity for improved efficiency.  Examples below will illustrate this.  Being able to isolate the interest parameter $\phi$ in the association is not always possible, it depends on the structure of the model, but there are strategies that may work even if the isolation is not strictly possible; see Section~\ref{SSS:bf}.

\subsection{Examples}
\label{SS:examples}

\subsubsection{Two normal means}
\label{SSS:twosample}

Consider two independent normal samples, i.e., 
\begin{equation}
\label{eq:baseline}
Y_{11},\ldots,Y_{1n_1} \iid \nm(\mu_1, \sigma^2) \quad \text{and} \quad Y_{21},\ldots,Y_{2n_2} \iid \nm(\mu_2, \sigma^2), 
\end{equation}
where the means, $\mu_1$ and $\mu_2$, and the sample sizes, $n_1$ and $n_2$, are possibly different, but the variance $\sigma^2$ is assumed to be the same in both groups; the more general case is considered in Section~\ref{SSS:bf}.  The full parameter in this case is $\theta=(\mu_1,\mu_2,\sigma)$, but the goal is only to determine if there is a difference between the two group means, a common marginal inference problem.  

As a first step, according to the general results in \citet{imcond}, reduce the dimension of the problem down via sufficient statistics.  Then the relevant association is 
\begin{align*}
\hat\mu_1 & = \mu_i + \sigma n_1^{-1/2} U_1, \\
\hat\mu_2 & = \mu_2 + \sigma n_2^{-1/2} U_2, \\
\hat\sigma^2 & = \sigma^2(n_1 + n_2 - 2)^{-1} U_3, 
\end{align*}
where $\hat\mu_1$ and $\hat\mu_2$ are the group sample means, $\hat\sigma^2$ is the usual pooled sample variance, $U_1$ and $U_2$ are independent standard normals, and $U_3$ is an independent chi-square with $n_1 + n_2 - 2$ degrees of freedom.  Interest is in $\phi = \mu_1 - \mu_2$ and manipulating the above association in an obvious way gives 
\[ \hat\phi = \phi + \hat\sigma \frac{n_1^{-1/2} U_1 + n_2^{-1/2} U_2}{\sqrt{(n_1 + n_2 - 2)^{-1} U_3}} \qquad \begin{array}{l}
\hat\mu_2 = \mu_2 + \sigma n_2^{-1/2} U_2, \\ \hat\sigma^2 = \sigma^2 (n_1 + n_2 - 2)^{-1} U_3, \end{array} \]
where $\hat\phi = \hat\mu_1 - \hat\mu_2$.  The latter two equations do not directly involve the interest parameter $\phi$, so if the goal is marginal inference on $\phi$, I have no need to predict $U_2$ and $U_3$, so these can be ignored.  And that particular feature of $(U_1,U_2,U_3)$ in the first equation has a known distribution, namely, a Student-t distribution with $n_1 + n_2 -2$ degrees of freedom.  So if $F$ denotes that distribution function, then the marginal association for $\phi$ can be written as 
\begin{equation}
\label{eq:marg.assoc}
\hat\phi = \phi + \hat\sigma (n_1^{-1} + n_2^{-1})^{1/2} \, F^{-1}(W), \quad W \sim \unif(0,1).
\end{equation}
If I now introduce the default random set $\S$ in \eqref{eq:default.prs}, which {\em is} optimal in this case, then I get a corresponding random set for $\phi$, 
\[ \Phi_x(\S) = \{ \varphi: w_y(\varphi) \in \S\}, \]
where $w_y(\varphi)$ is the solution of \eqref{eq:marg.assoc} for given data ``$y$'' and generic $\varphi$ value, i.e., 
\[ w_y(\varphi) = F\Bigl( \frac{\hat\phi - \varphi}{\hat\sigma \{n_1^{-1} + n_2^{-1}\}^{1/2}} \Bigr). \]
Then the plausibility is given by a probability with respect to the distribution $\prob_\S$ of the random set $\S$:
\begin{equation}
p_y(\{\varphi\}) = \prob_\S\{\Phi_y(\S) \ni \varphi\} = 1 - \bigl| 2 w_y(\varphi) - 1 \bigr|.
\end{equation}

\begin{figure}
\begin{center}
\subfigure[$\hat\phi=0.5$ and $\hat\sigma=1.25$.]{\scalebox{0.6}{\includegraphics{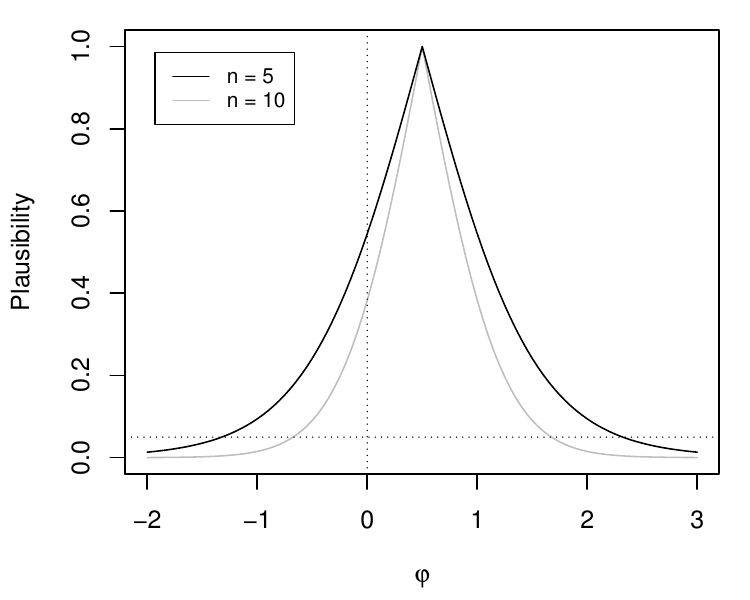}}}
\subfigure[$\hat\phi=1$ and $\hat\sigma=0.75$.]{\scalebox{0.6}{\includegraphics{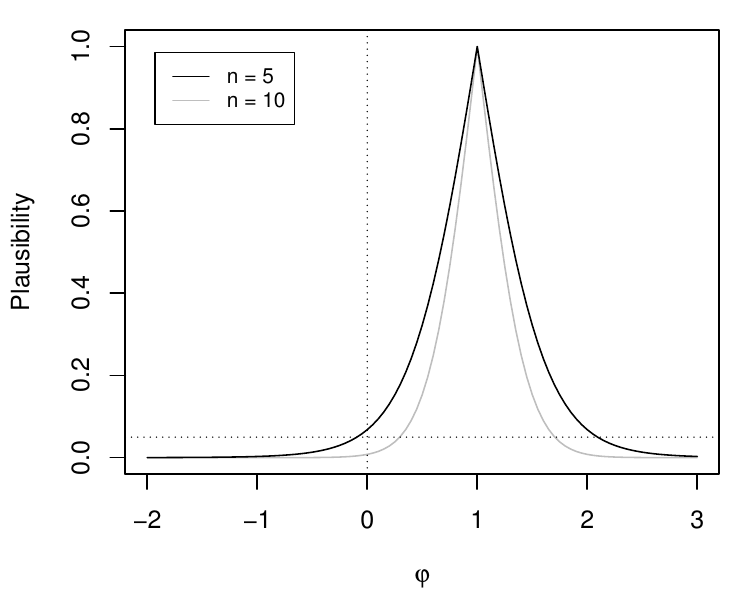}}}
\end{center}
\caption{Plots of the plausibility contour from Section~\ref{SSS:twosample} for the two cases and for the two sample size configurations in each case, all on the same scale.}
\label{fig:twosample}
\end{figure}

For illustration, I consider two $(\hat\phi, \hat\sigma)$ cases, one small signal-large noise and one large signal-small noise.  For simplicity of comparison only, I assume that $n_1=n_2=n$, and I consider $n=5$ and $n=10$ for each of the two cases above.  Figure~\ref{fig:twosample} shows the plausibility contour $p_y(\{\varphi\})$ for the two cases and two sample size configurations.  

First, in all four cases, the plausibility contour peaks at the estimate $\hat\phi$, and decays symmetrically as $\varphi$ moves farther from $\hat\phi$.  The interpretation is that the most likely value, the point estimator, is fully plausible---which is why it's used as an estimate---but values distant from the estimate should be less plausible.  Second, the more informative the data are, either because $\hat\sigma$ is smaller or because the sample size is larger, the more narrow the plausibility curve is.  This, again, agrees with the intuition that more informative data leads to more precise inference, i.e., a more narrow range of plausible values.  

What is especially interesting about the plausibility function plots in Figure~\ref{fig:twosample} is that the data analyst can immediately read off the answer to virtually any relevant question.  First, using the horizontal reference line at 0.05, one can easily identify the 95\% plausibility interval for $\phi$ as the set of all $\varphi$ with point plausibility at least 0.05.  Naturally, more informative data makes the interval narrower.  And the interpretation is clear: if a value, such as $\varphi=0$ is in the plausibility interval, then it is sufficiently plausible based on the data, otherwise it's not.  Note also that the plausibility function is on the appropriate scale for reading off these intervals; compare this to a Bayesian approach where the need to evaluate tail areas via integration makes it impossible to read off the credible interval from a plot of the posterior density.  Of course, there is nothing sacred about 0.05, one is free to set their own threshold for ``sufficiently implausible.''  

Second, since interest is in comparing the group means, a relevant hypothesis is $A=(-\infty, 0)$.  The reference line at $\varphi=0$ is helpful in this assessment.  In Panel~(a), for example, with $n=5$ (black line), I find that $p_y(\{0\})$ is between 0.5 and 0.6; the precise value is 0.545.  Therefore, since the point plausibility is increasing up to $\varphi=0$, I get $p_y(A) = 0.545$ and $b_y(A) = 0$, so no definitive judgment can be made about the hypothesis $A$ because belief is small but plausibility is not small.  Virtually the same conclusion would be reached in the slightly more informative case with $n=10$.  On the other hand, in Panel~(b), with $n=10$ (gray line), $p_y(\{0\})$ is effectively 0, so both $b_y(A)$ and $p_y(A)$ are very small, so I can infer $A^c$, i.e., the hypothesis $\{\phi < 0\} = \{\mu_1 < \mu_2\}$ is sufficiently implausible based on the data.  

\subsubsection{Normal coefficient of variation, cont.}

Reconsider the normal coefficient of variation example from Section~\ref{SS:cv}.  Following the general theory, take the minimal sufficient statistic $(\hat\mu, \hat\sigma^2)$, the sample mean and sample variance, respective, and build the association in terms of these:
\begin{align*}
\hat\mu & = \mu + \sigma n^{-1/2} U_1, \quad & U_1 & \sim \nm(0,1), \\
\hat\sigma^2 & = \sigma^2 U_2, \quad & U_2 & \sim (n-1)^{-1} \chisq(n-1),
\end{align*}
where $U_1$ and $U_2$ are independent.  The parameter of interest, $\phi = \sigma/\mu$, is a scalar but the auxiliary variable $(U_1,U_2)$ is two-dimensional, so we would like to further reduce the dimension of the latter.  An argument similar to that in Section~\ref{SSS:twosample} leads to a (marginal) association for the interest parameter $\phi$:
\[ n^{1/2} \hat\phi^{-1} = F_{n,\phi^{-1}}^{-1}(W), \quad W \sim \unif(0,1), \]
where $\hat\phi = \hat\sigma/\hat\mu$ and $F_{n,\psi}$ is the non-central Student-t distribution function, with $n-1$ degrees of freedom and non-centrality parameter $n^{1/2}\psi$.  Different random sets $\S$ for the unobserved $W$ may be considered (see below) but, for now, I'll consider that same ``default'' choice in \eqref{eq:default.prs}.  The calculation in Example~\ref{ex:gamma} leads to a plausibility contour  
\begin{equation}
\label{eq:cv.pl}
p_y(\{\varphi\}) = 1 - |2 F_{n,\varphi^{-1}}(n^{-1} \hat\phi^{-1}) - 1 |. 
\end{equation}
Since the random set is nested, the derived belief function is consonant and, hence, this plausibility contour fully determines the inferential model output.  

\begin{figure}
\begin{center}
\subfigure[$\mu=1$]{\scalebox{0.59}{\includegraphics{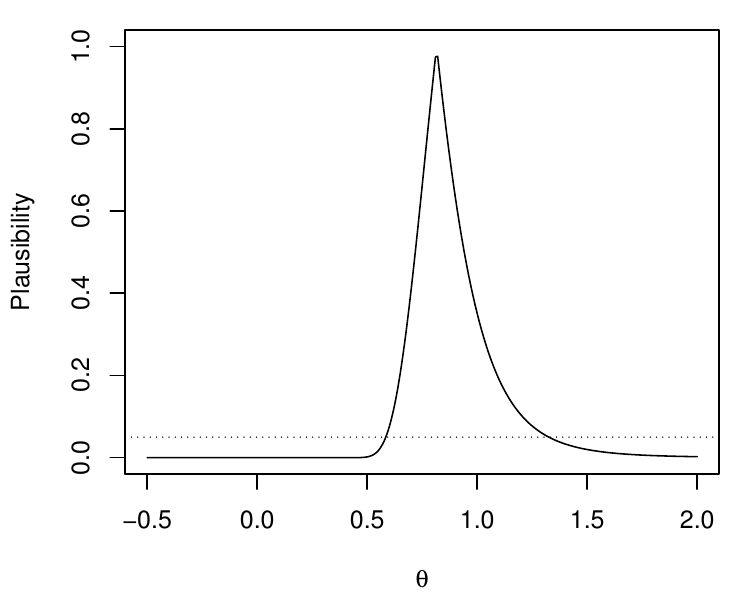}}} 
\subfigure[$\mu=0$]{\scalebox{0.59}{\includegraphics{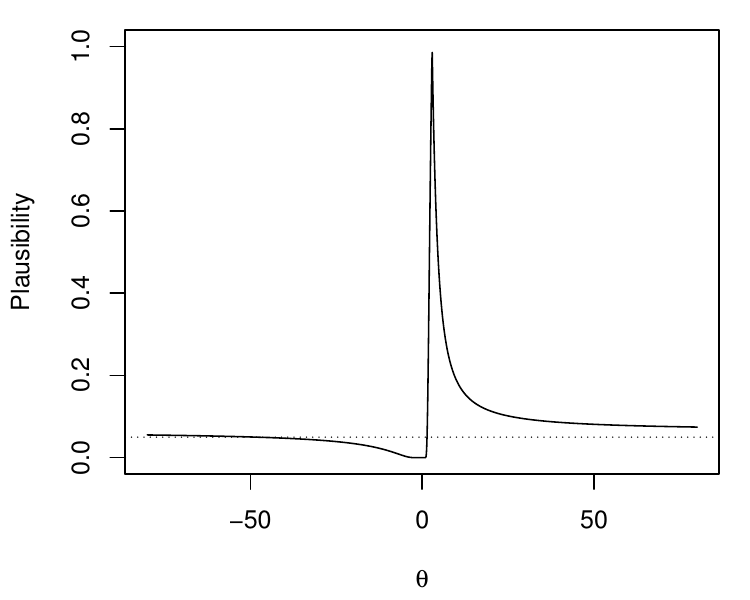}}}
\end{center}
\caption{Plots of the plausibility function for $\phi=\sigma/\mu$ based on samples of size $n=30$ from $\nm(\mu,\sigma^2)$ for $\sigma=1$ and two values of $\mu$.}
\label{fig:cv2}
\end{figure}

For given $\alpha \in (0,1)$, the $100(1-\alpha)$\% plausibility region for $\phi$ is given by 
\[ \{\varphi: p_y(\{\varphi\}) > \alpha\} = \{\varphi: \alpha / 2 < F_{n,\varphi^{-1}}(n^{1/2} \hat\phi^{-1}) < 1-\alpha/2 \}. \]
One can check directly that this plausibility region has coverage probability $1-\alpha$, but it follows as a consequence of Theorem~\ref{thm:freq}.   The Gleser--Hwang theorem applies to these plausibility regions, so having the nominal frequentist coverage implies that these will be unbounded with positive probability; see Figure~\ref{fig:cv2}, which shows two plausibility contours, one that vanishes in the tails, yielding bounded plausibility intervals, and one that does not.  Of course, unbounded regions are perfectly reasonable in this problem: given the inherent variability of the sample mean around $\mu$, if $\mu$ is close to zero, then arbitrarily large values of $\phi$ cannot be ruled out.  

\begin{figure}[t]
\begin{center}
\scalebox{0.7}{\includegraphics{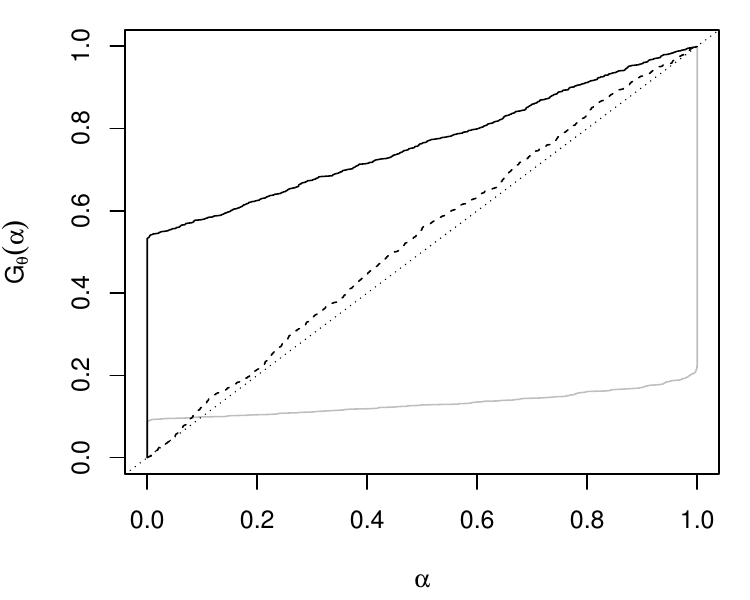}}
\end{center}
\caption{Plot of the distribution function $G_\theta(\alpha) = \prob_{Y|\theta}\{b_Y(A) \leq \alpha\}$ of the Bayesian posterior probability (gray), the valid inferential model (black, solid) with plausibility contour in \eqref{eq:cv.pl}, and the one based on the ``better'' random set (dashed), where $A=(-\infty,9]$ is a false hypothesis about the coefficient $\phi=\sigma/\mu$ when $\mu=0.1$ and $\sigma=1$.}
\label{fig:cv3}
\end{figure}

Next, I revisit the simulation from Section~\ref{SS:cv}.  Recall that $Y$ consists of $n=10$ iid observations from $\nm(\mu, \sigma^2)$, with $\mu=0.1$ and $\sigma=1$, so that $\phi=10$ is the true value.  Then I consider a false hypothesis $A=(-\infty, 9]$.  Here I want to compare the distribution of the belief function $b_Y(A)$ from above, as a function of $Y$, with that of the posterior distribution $\Pi_Y(A)$ from before.  Figure~\ref{fig:cv3} shows the same distribution function of $\Pi_Y(A)$ as in Figure~\ref{fig:cv1}, and also shows two other distribution functions: one for the belief function $b_Y(A)$ based on the inferential model construction laid out above; the other for the belief function based on the same formulation but with a ``better'' random set, namely, the one-sided interval $\S = [\tilde W, 1]$, with $\tilde W \sim \unif(0,1)$, suggested by Theorem~4 in \citet{imbasics}.  Both of the latter distribution functions stay above the diagonal line, because they're both valid, while the Bayes version falls below and, hence, is not valid.  That the one belief has distribution function closer to the diagonal line, without going below, is a sign of added efficiency; see Section~\ref{SS:efficiency}.  The takeaway message is that it's possible to construct valid (and even efficient) degrees of belief in this difficult problem, but non-additivity is essential.

\subsubsection{Behrens--Fisher problem}
\label{SSS:bf}

The so-called {\em Behrens--Fisher problem} is a fundamental one, a generalization of the problem in Section~\ref{SSS:twosample} where the two groups have their own unknown variances.  Despite its apparent simplicity, it turns out that the only ``exact'' procedures are those that depend on the order in which the data is processed, which is not fully satisfactory.  Standard solutions are given by \citet{hsu1938} and \citet{scheffe1970}; see, also, \citet{welch1938,welch1947}.  For a review, see \citet{kimcohen1998}, \citet{ghosh.kim.2001}, and \citet{fraser.wong.sun.2009}.  

Consider the model in \eqref{eq:baseline}, but where the two groups have their own unknown variance, $\sigma_1^2$ and $\sigma_2^2$, respectively.  Reduction via sufficiency leads to a baseline association
\begin{equation}
\label{eq:bf-aeqn1}
\hat\mu_k = \mu_k + \sigma_k \, n_k^{-1/2} \, U_{1,k}, \quad \text{and} \quad \hat\sigma_k = \sigma_k U_{2,k}, \quad k=1,2, 
\end{equation}
where $U_{1,k} \sim \nm(0,1)$ and $(n_k-1)U_{2,k}^2 \sim \chisq(n_k-1)$, independent, for $k=1,2$.  Define $f(\sigma_1,\sigma_2) = (\sigma_1^2/n_1 + \sigma_2^2/n_2)^{1/2}$ and note that \eqref{eq:bf-aeqn1} can be re-expressed as 
\begin{equation}
\label{eq:bf-aeqn2}
\hat\phi = \phi + f(\sigma_1,\sigma_2) \, V \quad \text{and} \quad \hat\sigma_k = \sigma_k U_{2,k}, \quad k=1,2, \quad V \sim \nm(0,1). 
\end{equation}
Plug in $\hat\sigma_k = \sigma_k U_{2,k}$ for $k=1,2$, to get
\begin{equation}
\label{eq:bf-aeqn3}
\frac{\hat\phi - \phi}{f(\hat\sigma_1,\hat\sigma_2)} = \frac{f(\sigma_1, \sigma_2)}{f(\sigma_1 U_{2,1}, \sigma_2 U_{2,2})} \, V, \quad \text{and} \quad \hat\sigma_k = \sigma_k U_{2,k}, \quad k=1,2. 
\end{equation}
If the coefficient on $V$ was free of $(\sigma_1,\sigma_2)$, it would be straightforward to marginalize to an association in terms of $\phi$ only, by ignoring the nuisance parameter components.  However, this simple approach is not possible in this challenging problem.  As an alternative, \citet{immarg} suggest a strategy in which the nuisance parameter components can be ignored but that the distribution of the auxiliary variable is stretched accordingly.  Skipping the details, they suggest a new dimension-reduced association for $\phi$,  
\[ \frac{\hat\phi-\phi}{f(\hat\sigma_1,\hat\sigma_2)} = G^{-1}(W), \quad W \sim \unif(0,1), \]
where $G$ is the distribution function of a Student-t distribution with $\max(n_1, n_2) - 1$ degrees of freedom.  Here, because of symmetry, the default random set $\S$ in \eqref{eq:default.prs} is optimal, and the plausibility contour can be evaluated exactly as in Example~\ref{ex:gamma} above.  Moreover, this generalized marginal inferential model is valid, which leads to an alternative proof of the conservative coverage properties of the Hsu--Scheff\'e confidence interval.  

\begin{figure}[t]
\begin{center}
\scalebox{0.70}{\includegraphics{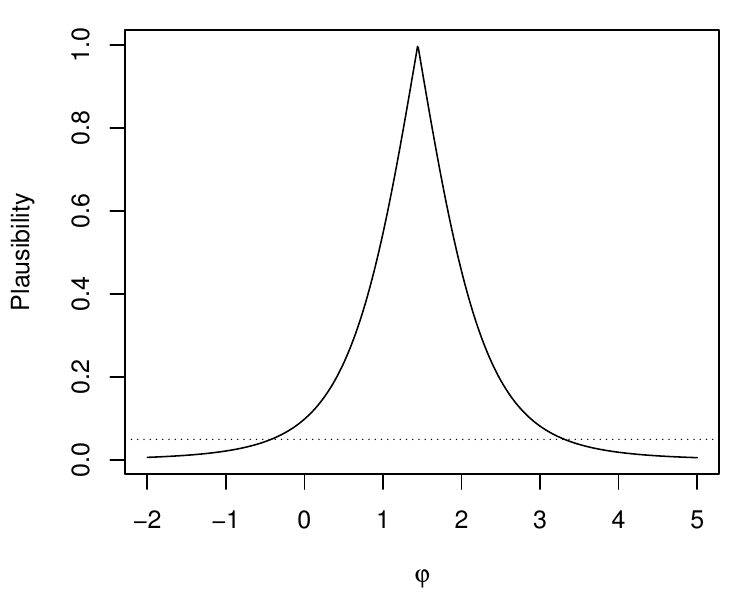}}
\end{center}
\caption{Plot of the plausibility contour for the example in Section~\ref{SSS:bf}.  Horizontal line at $\alpha=0.05$ determines the 95\% plausibility region for $\phi$.}
\label{fig:bf-fig}
\end{figure}

Data on travel times from home to work for two different routes are presented by \citet[p.~83]{lehmann1975}, summarized below:
\[ (n_1,\hat\mu_1,\hat\sigma_1^2) = (5, 7.58, 2.24) \quad \text{and} \quad (n_2, \hat\mu_2, \hat\sigma_2^2) = (11, 6.14, 0.07). \]
The goal is to investigate the difference between the mean travel times for the two routes.  Figure~\ref{fig:bf-fig} shows the plausibility contour and along with the horizontal cut at $\alpha=0.05$ that determines the 95\% plausibility interval for $\phi$.  As expected, the contour has a peak at the estimated mean difference, $\hat\phi = 1.44$, but the fact that $p_y(\{0\})$ exceeds 0.05 implies that the two routes having the same mean travel times is sufficiently plausible.  

\ifthenelse{1=1}{}{
\begin{table}[t]
\caption{Data on travel times for the example in Section~\ref{SSS:bf}.}
\label{table:bfdata}
\begin{center}
\begin{tabular}{c|ccc}
\hline
Route & $n$ & $\hat\mu$ & $\hat\sigma^2$ \\
\hline
1 & 5 & 7.58 & 2.24 \\
2 & 11 & 6.14 & 0.07 \\
\hline
\end{tabular}
\end{center}
\end{table}
}

\subsubsection{Two Poisson rates}
\label{ex:poisson}

Consider two independent Poisson samples, $Y_1 \sim \pois(n_1 \theta_1)$ and $Y_2 \sim \pois(n_2 \theta_2)$, where the ``sample sizes'' $n_1$ and $n_2$ are known but the rates $\theta=(\theta_1,\theta_2)$ are unknown.  The quantity of interest here is $\phi = \theta_1 / \theta_2$, the ratio.  Inference on $\phi$ has been explored in a number of papers, including \citet{krish.thomson.2004} and the references therein.  A similar problem arises in \citet{jin.li.jin.2015}, and I employed a conditioning strategy like the one below in \citet[][Sec.~5]{gim}.

The jumping off point is the realization that there are two data points---each with a corresponding auxiliary variables---but only a scalar interest parameter.  So reducing the dimension before building the inferential model would be advantageous.  One dimension-reduction strategy that can sometimes be employed is based on conditioning.  Here the approach is relatively straightforward and familiar.  That is, the conditional distribution of $Y_2$, given $Y_1 + Y_2 = k$, is binomial with parameters $k$ and 
\[ g(\phi) = \frac{n_2 \theta_2}{n_1 \theta_1 + n_2 \theta_2} = \{1 + (n_1/n_2) \phi\}^{-1}, \]
a one-to-one function of $\phi$.  So I just have to build an association in terms of the marginal distribution for $Y_1 + Y_2$ and also in terms of the aforementioned conditional.  Each will have an associated auxiliary variable but I will just ignore the one in the marginal part, focusing only on the association that describes the conditional distribution, which only involves the interest parameter.  And, fortunately, that analysis based on the conditional distribution is almost identical to that presented in Example~\ref{ex:binom1} above.  The only point of departure is the transformation between the binomial parameter $g(\phi)$ and $\phi$, but this is not a problem.  Here, however, I will also consider a random set other than the default considered in all the examples so far.  In particular, I will use the greedy random set construction laid out in \citet{impois} which, unfortunately, is too complicated to describe here; my reason for doing so is to illustrate what an {\em efficient} inferential model looks like in a discrete problem; see, also, Section~\ref{S:converse} below.

\begin{figure}
\begin{center}
{\scalebox{0.7}{\includegraphics{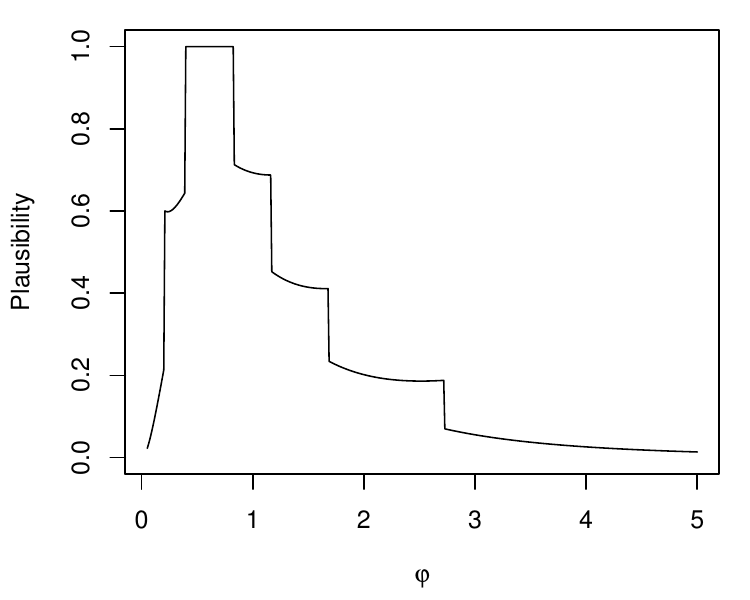}}}
\end{center}
\caption{Plot of the plausibility contour for the ratio $\phi=\theta_1/\theta_2$ of Poisson means based on the data as described in Example~\ref{ex:poisson}.}
\label{fig:poisson}
\end{figure}

As an illustration, suppose a builder purchases lumber from both Company 1 and Company 2.  Whenever the lumber is delivered, a sample is taken and inspected and, in the most recent deliveries, a sample of $n_1=30$ samples was taken from Company 1's lot and $y_1=2$ defectives were identified, while $n_2=35$ samples from Company 2's delivery yielded $y_2=4$ defectives.  Using the conditional-binomial formulation described above, with random set construction as in \citet{impois}, I get the plausibility contour for $\phi$ in Figure~\ref{fig:poisson}.  One might guess that these data can't distinguish between the two lumber suppliers in terms of their defective rates, and the plausibility contour plot confirms this since the plausibility at $\varphi=1$ is rather high.  What's especially interesting about the plot is that it's not smooth, there are both plateaus and corners.  My experience is that ``good'' random sets for discrete data problems will yield plausibility contours with this shape, but a clear understanding of this phenomenon is still lacking. 

\subsubsection{Fr\'echet shape}
\label{SSS:frechet}

Suppose that data $Y=(Y_1,\ldots,Y_n)$ are iid from a unit Fr\'echet distribution with shape parameter $\theta > 0$; usually there are location and scale parameters included but I'm leaving these out here only for simplicity.  The Fr\'echet distribution is {\em max-stable} and, therefore, is a canonical model in extreme value theory \citep[e.g.,][]{dehaan.ferreira.book}.  It's also a challenging problem because it lacks the nice exponential family structure that guarantees existence of a one-dimensional sufficient statistic.  So combining the different sources of information about $\theta$ in this case requires some non-standard techniques.  

To construct a valid inferential model, first recall that the distribution function for the Fr\'echet model is given by 
\[ F_\theta(y) = \exp(-y^{-\theta}), \quad y > 0, \]
and so this model can be described via an association
\begin{equation}
\label{eq:frec.assoc}
Y_i = (-\log U_i)^{-1/\theta}, \quad i=1,\ldots,n, 
\end{equation}
where $U_1,\ldots,U_n$ are iid $\unif(0,1)$ and the expression on the right-hand side above is just the solution $F_\theta^{-1}(u)$ of the equation $u=F_\theta(y)$.  As I hinted in Section~\ref{SS:dimred}, in such cases there will be features of the auxiliary variable that are fully observed, and that identifying those features---and suitably conditioning on them---is important for reducing the dimension and increasing the efficiency of the inferential model.  One might be able to guess directly what kind of features are observed (see below), but let me proceed with a formal construction which I think is informative.  

The jumping off point is the observation that, if $u_{y,\theta}$ is a solution to the system \eqref{eq:frec.assoc}, then those features, $\eta(u_{y,\theta})$, that don't vary with $\theta$ would correspond to observed features of $U$.  So then the goal would be to find a function $\eta$ that satisfies the following partial differential equation:
\[ \frac{\partial \eta(u_{y,\theta})}{\partial \theta} = 0. \]
Presently, it is unknown what problems will admit a differential equation that has a solution, but there are examples where solving this is possible, and there are other examples (see Section~\ref{SSS:cube} below) where some change-of-perspective is required in order to find a solution.  For the present case with $n$ iid samples, it makes sense to use the chain rule and re-express this differential equation as 
\begin{equation}
\label{eq:diffeq}
\frac{\partial \eta(u)}{\partial u} \Bigr|_{u=u_{y,\theta}} \cdot \frac{\partial u_{y,\theta}}{\partial \theta} = 0, 
\end{equation}
where $\partial u_{y,\theta} / \partial \theta$ is a $n$-vector that is easy to compute from the known form of $u_{y,\theta}$, and $\partial \eta(u) / \partial u$ is a matrix with $n$ columns and of rank $n-1$.  That is, $\eta$ is a map from the $n$-space of $u$'s to some $(n-1)$-dimensional subspace.  

Solving the above differential equation can proceed in various ways; techniques that I've found to be successful are guess-and-check and the {\em method of characteristics} \citep[e.g.,][]{polyanin2002} and here I'll illustrate the latter.  Write $m_\theta(u)$ for $\partial u_{y,\theta}/\partial \theta$ with $y=a(\theta,u)$ plugged in, so that it's only a function of $\theta$ and $u$.  In the present example, 
\[ m_\theta(u) = \theta^{-1} u \log u \log(-\log u), \quad u \in (0,1). \]
Then the method of characteristics suggests solving the following (nonparametric) system of ordinary differential equations:
\[ \frac{du_1}{m_\theta(u_1)} = \frac{du_2}{m_\theta(u_2)} = \cdots = \frac{du_n}{m_\theta(u_n)}. \]
The solution to such a system of equations looks like 
\[ z_\theta(u_1) = z_\theta(u_2) + c_1 = \cdots = z_\theta(u_n) + c_{n-1} \]
for constants $c_1,\ldots,c_{n-1}$, where $dz_\theta(u)/du = m_\theta(u)^{-1}$.  In the present example, 
\[ z_\theta(u) = \int m_\theta(u)^{-1} \,du = \theta \int \frac{1}{u \log u \log(-\log u)} \,du, \]
which, upon making a change-of-variable, $s=\log(-\log u)$, yields an anti-derivative 
\[ z_\theta(u) = \theta \log|\log(-\log u)|. \]
The constants $c_1,\ldots,c_{n-1}$ are intended to accommodate the rank constraint on the derivative of $\eta$, and suggests a sort of ``location constraint'' on the $z_\theta(u_i)$'s, so that the final solution is a function of the differences $z_\theta(u_i) - z_\theta(u_j)$.  That is, take 
\[ \eta(u) = M \, z_\theta(u), \]
where $z_\theta(u) = (z_\theta(u_1),\ldots,z_\theta(u_n))^\top$ and $M$ a matrix with constant vectors in its null space, so that the rows of $M$ are like contrasts.  Note that $\theta$ appearing as a multiplicative constant in $z_\theta(u)$ can be ignored here, so that $\eta(u)$ only depends on $u$.  From this final expression, it is easy to see that $\eta$ satisfies the above differential equation.  

In this particular example, the differences $z_\theta(u_i) - z_\theta(u_j)$ are proportional to the logarithm of the ratio $\log(-\log u_i)/\log(-\log u_j)$, which is easily seen to be equal to the ratio $\log y_i / \log y_j$.  Therefore, the feature $\eta(u)$ constructed above is observable, as intended.  The reader experienced with properties of the Fr\'echet distribution might have been able to anticipate this form of $\eta$ from the beginning, but I think this formal construction is still interesting and informative.  

\begin{figure}[t]
\begin{center}
\scalebox{0.7}{\includegraphics{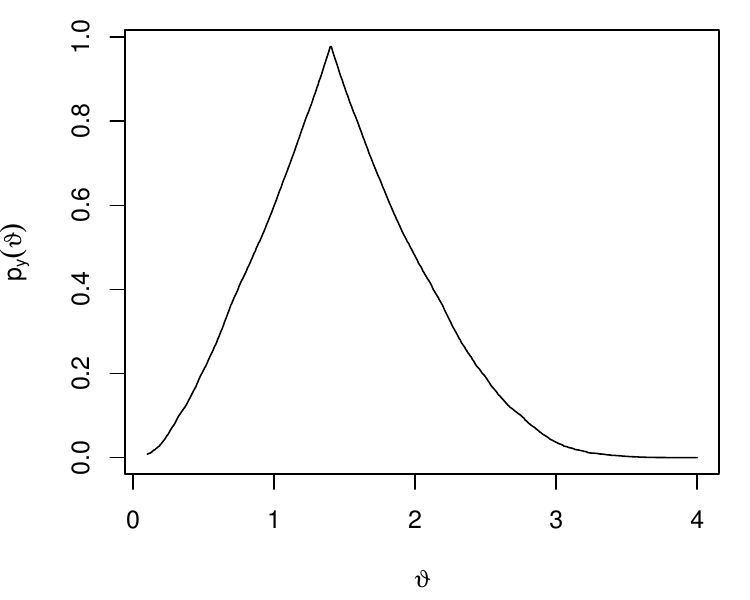}}
\end{center}
\caption{Plausibility contour for the Fr\'echet example; $n=2$, true $\theta=2$.}
\label{fig:frechet}
\end{figure}

To keep things simple here, let me consider only the $n=2$ case.  Then I can simplify the final result above and take 
\[ \eta(u) = \log(-\log u_2)/\log(-\log u_1). \]
I need another mapping, $\tau$, to complement $\eta$, and a suitable choice here is $\tau(u) = \log(-\log u_1)$.  This leads to a re-expressed ``dimension-reduced'' association 
\[ T(Y) = \theta^{-1} V_1 \quad \text{and} \quad H(Y) = V_2, \]
where $T(Y) = -\log Y_1$, $H(Y) = \log Y_2 / \log Y_1$, $V_1=\tau(U)$, and $V_2=\eta(U)$.  It is straightforward to find the conditional distribution of $V_1$, given $V_2=h$, where $h$ is the observed value of $H(Y)$, and then define a random set, $\S$, designed to predict draws from this conditional distribution.  That is, if $g_h$ is the density for the aforementioned conditional distribution, then a very reasonable---perhaps ``optimal''---random set is 
\[ \S = \{v_1: g_h(v_1) \geq g_h(\tilde V_1)\}, \quad \tilde V_1 \sim g_h. \]
With that the plausibility contour can be evaluated, via Monte Carlo, as 
\[ p_y(\{\vartheta\}) \approx \frac{1}{K} \sum_{k=1}^K 1\{ g_h(\tilde V_1^{(k)}) \geq g_h(\vartheta T(y))\}, \quad \tilde V_1^{(k)} \iid g_h. \]
A plot of this plausibility contour, based on a pair of samples ($n=2$) from $F_\theta$ with $\theta=2$, is shown in Figure~\ref{fig:frechet}. 

\subsubsection{Non-regular cube-root normal}
\label{SSS:cube}

Consider the following {\em cube-root normal} problem, 
\[ Y_1 = \theta + U_1 \quad \text{and} \quad Y_2 = \theta^{1/3} + U_2, \]
where $(U_1,U_2)$ consists fo two iid $\nm(0,1)$ random variables, and the goal is inference on $\theta$.  This example was considered in \citet{seidenfeld1992} for comparing fiducial and Bayes updates, and more recently in \citet{taraldsen.lindqvist.string}.  There's two observations carrying information about the same $\theta$, so there is an opportunity to reduce the dimension of the auxiliary variable to match that of the parameter.  No such reduction is provided by sufficiency in this non-regular example, but the differential equation strategy can be applied here to reduce the dimension.  

Write $u_{y,\theta}$ for the solution of the above system when $y$ and $\theta$ are fixed: 
\[ u_{y,\theta} = (y_1 - \theta, y_2 - \theta^{1/3})^\top. \]
As before, the goal is to find $\eta(u)$ to solve the differential equation
\[ \Bigl(\frac{\partial \eta(u)}{\partial u} \Bigr|_{u=u_{y,\theta}}\Bigr)^\top \frac{\partial u_{y,\theta}}{\partial \theta} = 0. \]
Unfortunately, it may not be possible to find a solution to this differential equation that doesn't depend on $\theta$.  But there's a way to get around this obstacle based on an idea of {\em localization} first presented in \citet{imcond}.  Let $\eta_\vartheta(u)$ be a function like before, depending on some generic parameter value $\vartheta$.  Localization requires that we solve the following differential equation:
\[ \frac{\partial \eta_\vartheta(u_{y,\theta})}{\partial \theta} = 0 \quad \text{when $\theta = \vartheta$}. \]
With this additional flexibility, the equation {\em can} be solved:
\[ \eta_\vartheta(u) = \vartheta^{2/3} u_2 - \tfrac13 u_1. \]
Now pair this $\eta_\vartheta$ with a function $\tau$ such that $u \mapsto (\tau(u), \eta_\vartheta(u))$ is one-to-one for each $\vartheta$.  Here I will take $\tau(u) = u_1$.  There is corresponding pair of functions $(T, H_\vartheta)$ on the $Y$-space, defined as 
\[ T(y) = y_1 \quad \text{and} \quad H_\vartheta(y) = \eta_\vartheta(u_{y,\vartheta}) = 3\vartheta^{2/3}(y_2 - \vartheta^{1/3}) - (y_1 - \vartheta). \]
Write out a ``dimension-reduced'' association
\[ T(Y) = \theta + \tau(U) \quad \text{and} \quad H_\vartheta(Y) = \eta_\vartheta(U). \]
The first equation depends on $\theta$ while the other does not.  This suggests writing $V_1 = \tau(U)$ and $V_2 = \eta_\vartheta(U)$ and using the conditional distribution of $V_1$, given the observed value $h_\vartheta = H_\vartheta(y_1,y_2)$ of $V_2$.  Of course, this conditional distribution is normal and can be worked out explicitly, although it's tedious; see below.  The symmetry of the normal auxiliary variable distribution suggests using the default random set in \eqref{eq:default.prs}. 

\begin{figure}[t]
\begin{center}
\scalebox{0.7}{\includegraphics{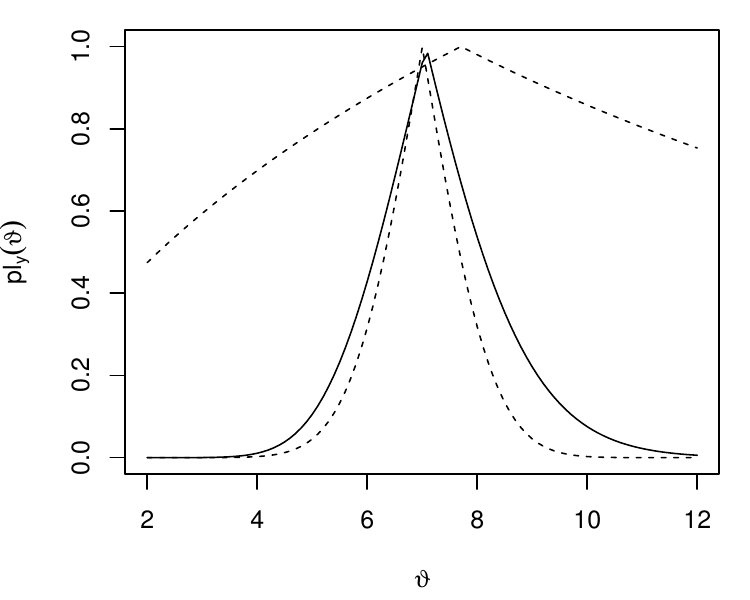}}
\end{center}
\caption{Plausibility contour for cube-root normal problem, true $\theta=7$.  Dashed lines are based on the individual data points, while the solid line is after combination using the differential equation strategy.}
\label{fig:cube}
\end{figure}

The above construction yields a different inferential model for each value $\vartheta$.  The key insight in \citet{imcond} was to ``glue'' these inferential models together on the plausibility contour scale.  That is, if $p_y(\cdot \mid \vartheta)$ was the plausibility contour for the $\vartheta$-dependent inferential model, then define the glued inferential model to have contour
\[ p_y(\{\vartheta\}) = p_y(\{\vartheta\} \mid \vartheta), \]
where the localization point varies with the argument to the plausibility contour.  Since the individual inferential models are valid, it is easy to show that the glued version is valid as well.  This gluing was used in other non-regular examples, including those in \citet{imunif} and in \citet{imvch}. In this example, with the default random set, the glued plausibility contour is 
\[ p_y(\{\vartheta\}) = 1-\Bigl|2 \Phi\Bigl( \frac{y_1 - \vartheta + \frac{1}{3\vartheta^{2/3}}(y_2 - \vartheta^{1/3})}{(1 + \frac{1}{9\vartheta^{4/3}})^{1/2}} \Bigr) - 1 \Bigr|. \]
For illustration, I simulate data from the above model with $\theta=7$.  Figure~\ref{fig:cube} shows a plot of this glued plausibility contour, along with those obtained by treating $y_1$ and $y_2$ separately.   Since $y_1$ is more informative about $\theta$ than $y_2$, it makes sense that the combined plausibility contour gives that observation more weight.  

\subsection{A sort of converse to Theorem~\ref{thm:freq}}
\label{S:converse}

So far, this section has focused on the construction of a valid inferential model, from which hypothesis tests and confidence regions with the desired frequentist error rate control can be immediately read off, as in Theorem~\ref{thm:freq}.  But, as is well known, one can get these decision procedures more directly, without full specification of degrees of belief.  In fact, a major advantage of a ``frequentist approach'' is the ability to break a problem down and jump directly to calibrated inferences on the quantity of interest.\footnote{I've recently been working on an approach by which one can construct marginal posteriors about parameters of interest directly, without first having a full posterior distribution.  We've been calling these {\em Gibbs posteriors} in \citet{syring.martin.image, syring.martin.scaling, syring.martin.mcid}.}  There's been lots of work over the years developing good---and in some cases optimal---tests and confidence regions so, for those researchers who prefer to focus on procedures with good properties, what does a valid inferential model have to offer?  

There are results available now which serve as a sort of converse to Theorem~\ref{thm:freq} above.  Roughly speaking, if there exists a test or confidence region that controls the frequentist error rate at any specified level $\alpha$, then there exists valid inferential model---potentially one that's {\em glued} in the sense of Section~\ref{SSS:cube}---whose plausibility function yields a test or confidence region at least as good (in terms of efficiency) as the given one.  The first such result was presented in \citet{impval} in the context of p-values and hypothesis testing.  More recently, I established a similar connection between valid inferential models and confidence regions by showing that, given a confidence region for a feature $\phi = \phi(\theta)$ of the full parameter, there exists a valid inferential model for $\theta$ whose corresponding marginal plausibility region for $\phi$ is at least as efficient as the given one \citep{imconf}.  Therefore, it's not just a superficial similarity between p-values/confidence regions and the valid belief/plausibility functions described here, there is a formal mathematical correspondence.  Moreover, this connection reveals that users of proper p-values and confidence regions have actually being working with the valid inferential model framework I described here all along, albeit unknowingly.  So, if this framework of valid inferential models can be considered a ``foundation,'' then this correspondence gives foundational support to the classical methods that are often considered to be ``ad hoc.''

\ifthenelse{1=1}{}{
My more recent results \citep{imconf} on the connection between valid inferential models and confidence regions go further.  I show that, given a confidence region for a feature $\phi = \phi(\theta)$ of the full parameter, there exists a valid inferential model for $\theta$ whose corresponding marginal plausibility region for $\phi$ is at least as efficient as the given one.  Moreover, there's an algorithm for carrying out this construction.  As a quick illustration, consider again the binomial problem in Example~\ref{ex:binom1}.  A classical procedure in the literature is the Clopper--Pearson confidence interval \citep[e.g.,][]{bcd2001}, and it is possible to construct a plausibility contour based on these intervals, which looks a lot like the plausibility function in Figure~\ref{fig:binom1}(a).  Using the framework in \citet{imconf}, I can convert the Clopper--Pearson intervals into a valid inferential model, and its plausibility contour is narrower, see Figure~\ref{fig:binom2}, hence the derived plausibility intervals are more efficient. 

\begin{figure}[t]
\begin{center}
\scalebox{0.7}{\includegraphics{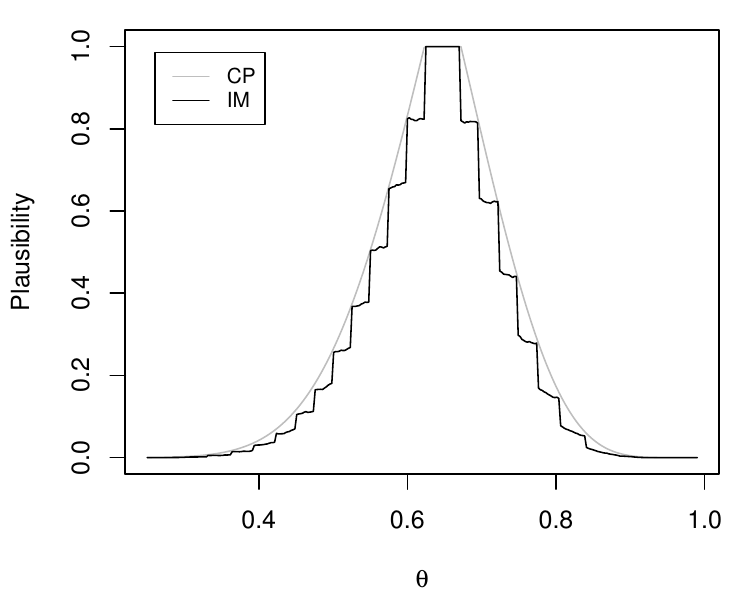}}
\end{center}
\caption{Plausibility contours corresponding to the Clopper--Pearson confidence interval for the binomial problem discussed in Section~\ref{S:converse}, along with that derived from the valid inferential model construction in \citet{imconf}; $n=20$ and $y=13$.}
\label{fig:binom2}
\end{figure}

From Section~\ref{S:intro}, calibration is essential for reliable statistical inference.  This reliability can be achieved by working with a valid inferential model or, to a certain extent, by developing tests and confidence intervals with frequentist error rate control for each individual task.  The take-away message here is that there's virtually no loss of flexibility or efficiency in taking the former perspective since it is always possible to construct an inferential model whose tests and confidence regions are at least as efficient.  
}

\section{More remarks}
\label{S:remarks2}

\subsection{A new take on the fiducial argument}
\label{SS:fiducial}

In Section~\ref{SS:bff} I described Fisher's version of the fiducial argument, circa 1930.  That argument is delightful, but requires too much structure.  \citet{dempster1963, dempster1964} describes it in a different way, which I'll summarize here and then relate back to what's going on the valid inferential model construction.  

Consider the baseline association \eqref{eq:assoc} from which data can be generated, i.e., sample $U \sim \prob_U$, plug into the mapping, for a given $\theta$, and return $Y=a(\theta, U)$.  Clearly, $U$ and $Y$ are dependent; in fact, they're very tightly linked, via a deterministic mapping, so there's no joint distribution as usual.  If one had a ``conditional distribution of $U$, given $Y=y$,'' then one could just sample, say $\tilde U$, from that distribution, plug it into $y = a(\theta, \tilde U)$, and solve for $\theta$ to get a draw from its ``posterior distribution.''  The problem, however, is that there's no such conditional distribution, at least not one that's useful.  That is, if the association mapping $a$ is sufficiently nice, then observing $Y=y$ immediately fixes $U$ at $u_{y,\theta}$, the solution to the equation $y=a(\theta,u)$, where $\theta$ is the true value.  In other words, $\prob_{U|Y=y} = \delta_{u_{y,\theta}}$, where $\delta_x$ denotes a point mass distribution at $x$.  Of course, this ``conditional distribution'' is useless because it depends on the true-but-unknown $\theta$.   

If not from probability calculus, where then does the fiducial distribution come from?  The word ``fiducial'' means {\em based on trust}, and Dempster made Fisher's leap of faith clear.  Indeed, Fisher just defines---or assumes---a particular conditional distribution, namely, 
\[ \prob_{U|Y=y} \overset{\text{\tiny set}}{\, = \,} \prob_U. \]
As Dempster described it, Fisher's suggestion was to ``continue to regard'' $U$ as if its distribution were $\prob_U$ even after observing $Y=y$, effectively ignoring the strong dependence between them.  Fisher's ``continue to regard'' idea is clever, and has been used for many years by those working in the fiducial domain.  

The point here so far is that, no surprise, it's not possible to construct a genuine posterior/conditional distribution out of thin air.  And if the fiducial distribution isn't built up according to the rules of probability, then it's again no surprise that it can break down when one puts the pressure of the probability calculus on it.  

At a high level, what causes problems in Fisher's argument is the fact that $Y$ and $U$ are too tightly tied together for ordinary probability calculus to handle.  That is, if $U$ and $Y$ are deterministically linked, as in \eqref{eq:assoc}, in a $\theta$-dependent way, then conditioning won't work.  The question, then, is how to break that tight link?   Dempster's efforts, in a series of classic papers \citep{dempster1966, dempster1967, dempster1968a, dempster1968b, dempster1969}, to develop what he now calls a {\em state space model} that describes all the unknowns and the various manipulations (e.g., Dempster's rule of combination) are carried out is, I think, in this vein.  

The approach described in Section~\ref{S:valid.im} has similar motivation and I want to focus the rest of this subsection on that.  As described above, the difficulties with fiducial inference arise because there's too tight of a link between $Y$ and $U$.  On the other hand, our original inspiration for the above inferential model construction was this tight connection.  That is, if identification of $u_{y,\theta}$ were possible, then the true $\theta$ could be similarly identified.  So the tight link between $Y$ and $U$ makes it possible to convert the problem of inference about $\theta$ to one of predicting or guess $u_{y,\theta}$.  This new task might not seem any easier than the previous, but the difference is that there's information available, in the generator $\prob_U$, about roughly where $u_{y,\theta}$ is located.  That information is used to generate a ``net''---a random set---to be cast in order to catch the fixed $u_{y,\theta}$.  There's no longer any concern about the link between $Y$ and $U$ because these random sets live on an entirely different probability space, independent of data and original auxiliary variable.  Since the target $u_{y,\theta}$ is a draw from $\prob_U$, it makes sense to design these random sets that they can hit targets drawn from $\prob_U$ with some amount of frequency.  If that can be arranged, then the set of solutions $\vartheta$ to $y=a(\vartheta,u)$ corresponding to $u$ in that set will inherit exactly the same frequency of containing the true $\theta$, and that's basically the proof of Theorem~\ref{thm:valid}.  Of course, I would also have the option to choose the random set to be a singleton filled with a draw from $\prob_U$, like Fisher's and Dempster's solutions, but this typically will not satisfy the desired validity property.  See \citet{liu.martin.wire}.  

Therefore, the construction presented in Section~\ref{S:valid.im} above is very much in line with Fisher's fiducial argument, most importantly, because the goal is to construct degrees of belief about the parameter where none had existed before, at least not formally.  The key difference, however, is that the inferential model construction doesn't simply settle with whatever comes out of the ``continue to regard'' step, it proactively seeks out a random set formulation with the desired validity property.

\subsection{Efficiency and optimality}
\label{SS:efficiency}

In the above examples, I have primarily made use of the default random set $\S$ defined in \eqref{eq:default.prs}.  In some examples, like the gamma model in Example~\ref{ex:gamma}, this choice was for simplicity, but in others I claimed that the default $\S$ was ``optimal'' in some sense.  Here I want to say a few words about {\em efficiency} in the inferential model context, and the corresponding notion of {\em optimality}, which means ``most efficient.''  

In classical statistics, efficiency can be expressed in terms of small (asymptotic) variance, e.g., between two confidence intervals that both attain the nominal coverage probability, the narrower one is more efficient.  Recall that the size of the plausibility region \eqref{eq:plausibility.region} is controlled by the spread of plausibility contour which, in turn, is controlled by the random set $\S$.  So, in the inferential model context, efficiency will be characterized or determined by the random set $\S$ or, rather, its distribution $\prob_\S$.  

For the moment, I will write $b_y(\cdot ; \S)$ to indicate that the random set $\S$ was used in the construction of the inferential model.  Then I'll say that $\S$ is at least as efficient as $\S'$ at hypothesis $A$, subject to both being valid, if $b_y(A; \S) \geq b_y(A; \S')$ for all $y$.\footnote{Asking for the inequality to hold uniformly in $y$ is indeed a strong type of efficiency, perhaps too strong.  Replacing this inequality with a distributional inequality is possible, and probably better.  The one result quoted here, however, holds under this strong condition.}  This implies that $b_Y(A; \S)$ is stochastically no smaller than $b_Y(A; \S')$ or, in other words, the spread of $p_Y(\vartheta; \S)$ would tend to be narrower than that of $p_Y(\vartheta; \S')$.  

A first basic result, presented in \citet{imbasics}, is a complete-class theorem which says that, for any random set $\S'$, there exists a {\em nested} random set $\S$ that is at least as efficient.  This explains why the inferential output is always a possibility function and not a general belief function. A nested random set has a fixed center, or core, so only its size and shape can vary.  But once the overall shape is fixed, the size is determined by the need to meet the validity condition.  To see this, consider, like in the discussion around \eqref{eq:h.prs}, a function $h: \UU \to \RR$ and a random set defined by 
\[ \S = \{u \in \UU: h(u) \leq s(\tilde U)\}, \quad \tilde U \sim \prob_U, \]
where $s: \UU \to \RR$ is another function that controls the size of $\S$.  If $s$ isn't at least as large as $h$ everywhere, then validity might fail, and if $s > h$ uniformly, then $\S$ is too large, so the best choice is to take $s=h$.  It is in this sense that the size of the random set is determined by its shape.  So the choice of an efficient $\S$ boils down to choosing the appropriate shape.  In certain one-dimensional situations, the optimal shape is known \citep{imbasics}.  For example, in symmetric location problems, the optimal random set is itself symmetric, hence \eqref{eq:default.prs} is optimal when $U$ is $\unif(0,1)$.  

Efforts to understand what makes a ``good'' random set are ongoing, but my current intuition is that the support of $\S$, i.e., the collection of possible realizations, should match the contours of the distribution $\prob_U$, like in Section~\ref{SSS:frechet}.  What makes formalization of this intuition non-trivial is that the association mapping $a: \Theta \times \UU \to \YY$ plays a role, as does the particular hypothesis $A \subset \Theta$ that one is interested in.  But I believe that the optimal random set's dependence on the association mapping $a$ is actually beneficial in the sense that it will all or mostly eliminate the inferential model's dependence on the functional form of the association, so that effectively no extra ``$\ldots$'' inputs beyond the data and statistical model would be needed for its construction.

\subsection{Practicalities}
\label{SS:practicalities}

There are certainly some differences---philosophical or otherwise---between what has been presented above and the more traditional statistical analyses, but the practical differences are not as wide as one might think.  Indeed, in those classical examples, such as linear regression, analysis of variance, etc., the recommended inferential model will have plausibility contour that yields hypothesis tests and confidence intervals that agree with those presented in the standard textbooks.  So, very little in terms of day-to-day operations of statistics requires a fundamental change.

From an educational point of view, at least for those courses intended for non-statistics majors in social or biological sciences, instructors could decide if introducing the plausibility contour used to produce the standard hypothesis tests and confidence intervals is appropriate for their students.  But I firmly believe that there is an advantage to using language like {\em plausibility} when communicating statistical ideas to these students compared to the more common {\em probability}, {\em confidence}, etc.  The reason is that these students know just enough probability to be dangerous.  For example, if a student confuses p-value for a probability, as many do, then he is likely to incorrectly interpret a small p-value as direct support for the truthfulness of the alternative hypothesis when, in reality, small p-value is only indirect support for the alternative, i.e., plausibility of the alternative is large but its belief could be small.  That student may not be able to understand the sub-additivity properties of plausibility functions, but if he knows that a p-value represents plausibility instead of probability, and that the former is a different object that can't be manipulated by the same rules as the latter, then he's likely not to make the aforementioned mistake.  And in light of the results summarized in Section~\ref{S:converse}, ``confidence'' is (mathematically) in direct correspondence to plausibility, so there is no need to use the former word if the latter is both satisfactory and sufficiently general.  For a different set of students, with more of a mathematical background, ideas about how to carry out the entire intermediate-level statistics theory course based on plausibility contours are presented in \citet{pvalue.course}.  

The majority of real-life data analyses consist of applying those standard methods taught in the course described above, so little or nothing would change for practicing statisticians.  And it would be straightforward to modify the existing statistical software programs to return a plot of the plausibility contour if users request it.  

That (certain features of) the recommended inferential model agree with traditional analyses in some problems, of course, does not imply that the solutions would agree in all problems.  Indeed, there are already a number of examples where the inferential model constructed as above produce different and usually better solutions than other methods.  More on the potential of valid inferential models is presented in Section~\ref{S:discuss}.

\section{Conclusion}
\label{S:discuss}

In this paper I have presented a summary of my current views on statistical inference.  The key foundational point is that, if one agrees with Reid and Cox that ``calibrated inferences seem essential,'' then additive degrees of belief are not appropriate for describing uncertainty.  That is, additive beliefs are afflicted with false confidence and, therefore, at least for some hypotheses, Theorem~\ref{thm:fct} says that there is a systematic and practically problematic bias.  Two examples in Section~\ref{S:additive} showed that the hypotheses afflicted need not be trivial ones either, so there is a genuine risk.  For those data analysts who are reluctant to move away from their familiar additive beliefs, it's important to be able to identify which hypotheses are at greatest risk of false confidence so they can mitigate their damages.  Towards this, those ``impossible'' problems highlighted in Section~\ref{SS:cv} should hold some insights.  But, as I explained in Section~\ref{S:converse}, my suggested move to non-additive inferential models does not necessarily affect everyday statistical practice; in fact, the results summarized there that characterize suitable p-values and confidence intervals as derivatives from a valid inferential model give foundational support to those classical methods.  

My focus here has been on the situation where no prior information is available, which might be unrealistic.  Is {\em nothing} really known about $\theta$?  What's typically done is either the available information is ignored, because it falls short of being a full prior distribution, or the inevitable information gaps are filled in by something artificial, e.g., default priors.  An ideal situation is that one can incorporate {\em exactly} the available information, no omissions and no additions.  Fortunately, there are ways in which this information can be incorporated into the inferential model framework.  Recently, \citet{chuanhai.partial} developed a valid inferential model for the case where prior information is available about nuisance parameters, in the form of a genuine Bayesian prior distribution, but not about the parameter of interest.  \citet{imexpert} presents preliminary results on a different angle, one where potentially incomplete prior information is available and to be used but, through the combination mechanism, will only lead to an improved inferential model, i.e., validity will be maintained and no efficiency will be lost, regardless of whether the prior is ``right'' or ``wrong,'' but efficiency will be gained when the prior is ``right.''   

Another interesting direction I'm currently pursuing, one that has a lot in common with the partial prior problem described above, is that of {\em model uncertainty}.  For example, suppose I am entertaining a collection of possible models, indexed by $M \in \M$, each having a model parameter $\theta_M$.  A Bayesian framework is equipped to handle such a problem, at least in principle, but difficulties quickly arise.  For example, one might have some vague prior information about the model itself, e.g., a quantification of Occam's razor which puts more weight on simpler models, but rarely does one have genuine prior information about the individual $\theta_M$'s, and the usual rules about non-informative priors don't apply when multiple models are under consideration.  An idea along these lines is to carry out the marginalization steps as described in Section~\ref{SS:dimred} above to isolate $M$ and, if available, introduce a sort of ``prior'' for $M$ along the lines described in the previous paragraph.  Preliminary work in this direction can be found in \citet{martin.model.isipta}.  

That the inferential model constructed in Section~\ref{S:valid.im} is able to achieve the validity property is exciting, but since there's no free lunch, there must also be some limitations.  At present, there are still many questions about this framework that I don't know how to answer, but I think this is more of an opportunity than a limitation.  As far as limitations are concerned, I'll mention two here.  First, since validity and efficiency are at odds with one another, and since validity is inflexible, some loss of efficiency is possible.  That is, if I construct a generic inferential model about $\theta$, then it's possible that the inference based on $b_y(A)$ or $p_y(A)$ for some specific $A$ is not as efficient as that based on some other approach or method that targets $A$ specifically.  If, however, the inferential model is tailored to the specific hypothesis (or hypotheses) of interest, then my experience is that efficiency can be maintained.  That leads to the second limitation, namely, difficulty.  Presently, it is non-trivial in some cases to carry out the necessary dimension reduction steps to achieve both validity and efficiency.  I will continue to think about ways to automate the process in certain ways, it's just hard to imagine what that ``black box'' would look like.  The Bayesian approach faced similar challenges, and its present popularity is entirely due to the advent of Monte Carlo methodology that virtually automated the computations that had previously been impossible.  What will the next breakthrough be?  Is a similar automation of inferential model construction and belief function computation possible?

A question that has come up several times in conversations with researchers about the material presented herein is if there is a version of {\em decision theory} under the inferential model umbrella.  There have, of course, been considerable efforts over the years to build up a belief function parallel to the decision theory of \citet{savage1972} based on probabilities; see the recent survey by \citet{denoeux.decision.2019}.  These can be applied off the shelf to the belief functions resulting from the inferential model construction described in Section~\ref{S:valid.im} here, but it is not clear---at least to me---what role validity plays in that context.  Intuitively, the expected utility might be affected by false confidence lurking in the posterior distribution, so there could be some advantage to working with a valid belief function that is immune to false confidence, but more work is needed.

Now is an exciting time for statisticians and data scientists: lots of job opportunities in both industry and academia, interesting (``big-data'') research problems with important real-life consequences, and a new generation of STEM students eager to be trained.  But there are also dangers.  As problems become bigger and more complex, we may reach a point where our intuition is no longer reliable.  For example, dependent data is commonplace in certain applications, and dealing with that dependence is by now familiar.  But modern network problems bring an entirely different meaning to ``dependent data,'' i.e., {\em the data is the dependence}.  Can we really expect our intuitions, developed largely based on asymptotic approximations in relatively simple cases, to safely guide us through these new complex problems?  After all, false confidence is lurking in this problem too, and it's surely even more difficult to spot amidst all the complexity.  

I'd like to end, just like I began, by quoting from \citet{barnard.book.1985}:
\begin{quote}
{\em It seems to be useful for statisticians generally to engage in retrospection at this time, because there seems now to exist an opportunity for a convergence of view on the central core of our subject.}
\end{quote}
Before having put together the results presented in \citet{balch.martin.ferson.2017} and the valid inferential model construction in \citet{imbook}, I would've said that this ``convergence of view'' ship had sailed.  The difference is that now we understand what can go wrong in practice with additivity, not just philosophies and generalities, along with a normative approach for the construction of valid inferential models. 

\section*{Acknowledgments}

I am sincerely grateful to Chuanhai Liu for his contributions to an earlier version of this manuscript and, more generally, for sharing so many insights and experiences over the years.  I have also benefited greatly from conversations about probability and statistics with Michael Balch and Harry Crane, and I thank two Statistics PhD students at North Carolina State University, namely, Mr.~Jesse Clifton and Mr.~Alex Cloud, for challenging me with tough questions.  Last but not least, I thank Thierry Denoeux for the invitation to give a lecture at the {\em BELIEF/SMPS 2018} conference in Compi\'egne, the Guest Editors for the invitation to contribute to the corresponding special issue of the {\em International Journal of Approximate Reasoning}, and the two anonymous reviewers for their efforts to read and provide helpful feedback on this lengthy manuscript.

\bibliographystyle{apalike}
\bibliography{/Users/rgmarti3/Dropbox/Research/mybib}

\end{document}